\pdfoutput=1
\RequirePackage{fix-cm}
\documentclass[smallextended,envcountsame,envcountsect,numbook]{svjour3}       % onecolumn (second format)

\smartqed  % flush right qed marks, e.g. at end of proof
\usepackage{graphicx}
\usepackage[T1]{fontenc}

\usepackage{amssymb,amsmath,amscd}

\usepackage{mathrsfs}
\usepackage{amsthm}
\usepackage{latexsym}
\usepackage{enumerate}

\usepackage{stmaryrd}

\usepackage{graphicx}
\usepackage{tikz-cd}

\providecommand{\constsh}[1]{\underline{#1}}

\newcommand{\style}[1]{{\emph{#1}}}  %\index{{#1}}}    % DEFINITION STYLE

\newcommand{\ie}{{\em i.e.}}
\newcommand{\eg}{{\em e.g.}}

\newcommand{\face}{{\trianglelefteqslant}}%\lhd}}  % FACE RELATION

\newcommand{\Field}{{\Bbbk}}
\newcommand{\Lift}{{\varphi}}

\newcommand{\Vspace}{{\mathbb V}}
\newcommand{\Thresh}{{\kappa}}

\newcommand{\R}{\mathbb{R}}
\newcommand{\Pb}{\mathbb{P}}

\newcommand{\C}{\mathbb{C}}

\newcommand{\E}{\mathbb{E}}

\newcommand{\Fc}{\mathcal{F}}
\newcommand{\Gc}{\mathcal{G}}

\DeclareMathOperator{\Hom}{Hom}

\DeclareMathOperator{\img}{im}

\DeclareMathOperator{\id}{id}
\DeclareMathOperator{\tr}{tr}
\DeclareMathOperator{\rank}{rank}
\DeclareMathOperator{\im}{im}

\DeclareMathOperator{\codim}{codim}

\DeclareMathOperator{\st}{st}

\DeclareMathOperator{\Vect}{\mathbf{Vect}}
\DeclareMathOperator{\Hilb}{\mathbf{Hilb}}

\providecommand{\abs}[1]{\left\lvert#1\right\rvert}
\providecommand{\norm}[1]{\lVert#1\rVert}

\providecommand{\ip}[1]{\langle #1 \rangle}

\providecommand{\abs}[1]{\left\lvert#1\right\rvert}
\providecommand{\norm}[1]{\lVert#1\rVert}

\providecommand{\ip}[1]{\langle #1 \rangle}

\def\@lbibitem[#1]#2{\item[{[#1]}\hfill]\if@filesw
      {\let\protect\noexpand
       \immediate
       \write\@auxout{\string\bibcite{#2}{#1}}}\fi\ignorespaces}

\newcommand{\etalchar}[1]{$^{#1}$}

\journalname{}%Journal of Applied and Computational Topology}
\begin{document}

\title{Toward a Spectral Theory of Cellular Sheaves
}

\author{Jakob Hansen
         \and
        Robert Ghrist
}

\institute{J. Hansen \at
             Department of Mathematics \\
             University of Pennsylvania \\
             David Rittenhouse Lab. \\
209 South 33rd Street \\
Philadelphia, PA 19104-6395 \\
              \email{jhansen@math.upenn.edu}
           \and
           R. Ghrist \at
             Departments of Mathematics and Electrical \& Systems Engineering \\
             University of Pennsylvania \\
             David Rittenhouse Lab. \\
209 South 33rd Street \\
Philadelphia, PA 19104-6395 \\
              \email{ghrist@math.upenn.edu}
}

\date{Received: date / Accepted: date}

\maketitle

\begin{abstract}
This paper outlines a program in what one might call {\em spectral sheaf
theory\/} --- an extension of spectral graph theory to cellular sheaves. By
lifting the combinatorial graph Laplacian to the Hodge Laplacian on a cellular
sheaf of vector spaces over a regular cell complex, one can relate spectral
data to the sheaf cohomology and cell structure in a manner reminiscent of
spectral graph theory. This work gives an exploratory introduction, and
includes discussion of eigenvalue interlacing, sparsification, effective
resistance, synchronization, and sheaf approximation. These results and
subsequent applications are prefaced by an introduction to cellular sheaves and
Laplacians.
\end{abstract}

\pagebreak
% ==========================================================================================
\section{Introduction}\label{sec:intro}
% ==========================================================================================
In spectral graph theory, one associates to a combinatorial graph additional
algebraic structures in the form of square matrices whose spectral data is then
investigated and related to the graph. These matrices come in several variants,
most particularly degree and adjacency matrices, Laplacian matrices, and
weighted or normalized versions thereof. In most cases, the size of the
implicated matrix is based on the vertex set, while the structure of the matrix
encodes data carried by the edges.

To say that spectral graph theory is useful is an understatement. Spectral
methods are key in such disparate fields as data
analysis~\cite{belkin_laplacian_2003,coifman_diffusion_2006}, theoretical
computer science~\cite{hoory_expander_2006,cvetcovic_graph_2011}, probability
theory~\cite{lyons_probability_2016}, control
theory~\cite{bullo_lectures_2018}, numerical linear
algebra~\cite{spielman_nearly_2014}, coding
theory~\cite{spielman_linear-time_1996}, and graph theory
itself~\cite{chung_spectral_1992,brouwer_spectra_2012}.

Much of spectral graph theory focuses on the Laplacian, leveraging its unique
combination of analytic, geometric, and probabilistic interpretations in the
discrete setting. This is not the complete story. Many of the most well-known
and well-used results on the spectrum of the graph Laplacian return features
that are neither exclusively geometric nor even combinatorial in nature, but
rather more qualitative. For example, it is among the first facts of spectral
graph theory that the multiplicity of the zero eigenvalue of the graph
Laplacian enumerates connected components of the graph, and the relative size
of the smallest nonzero eigenvalue in a connected graph is a measure of
approximate dis-connectivity. Such features are topological.
%, and it is here that the present inquiry begins. TOO FLOWERY

There is another branch of mathematics in which Laplacians hold sway: Hodge
theory. This is the slice of algebraic and differential geometry that uses Laplacians on
(complex) Riemannian manifolds to characterize global features. The classical
initial result is that one recovers the cohomology of the manifold as the
kernel of the Laplacian on differential forms~\cite{abraham_manifolds_1988}.
For example, the dimension of the kernel of the Laplacian on $0$-forms
($\mathbb{R}$-valued functions) is equal to the rank of $H^0$, the $0$-th
cohomology group (with coefficients in $\mathbb{R}$), whose dimension is the
number of connected components. In spirit, then, Hodge theory categorifies
elements of spectral graph theory.

Hodge theory, like much of algebraic topology, survives the discretization from
Riemannian manifolds to (weighted) cell
complexes~\cite{eckmann_harmonische_1945,friedman_computing_1998}. The
classical boundary operator for a cell complex and its formal adjoint combine
to yield a generalization of the graph Laplacian which, like the Laplacian of
Hodge theory, acts on higher dimensional objects (cellular cochains, as opposed
to differential forms). The kernel of this discrete Laplacian is isomorphic to
the cellular cohomology of the complex with coefficients in the reals,
generalizing the connectivity detection of the graph Laplacian in grading zero.
As such, the spectral theory of the discrete Laplacian offers a geometric
perspective on algebraic-topological features of higher-dimensional complexes.
Laplacians of higher-dimensional complexes have been the subject of recent
investigation~\cite{parzanchevski_high_2013,steenbergen_towards_2013,horak_spectra_2013}.

This is not the end. Our aim is a generalization of both spectral graph theory
and discrete Hodge theory which ties in to recent developments in topological
data analysis. The past two decades have witnessed a burst of activity in
computing the homology of cell complexes (and sequences thereof) to extract
robust global features, leading to the development of specialized tools, such
as persistent homology, barcodes, and more, as descriptors for cell
complexes~\cite{carlsson_shape_2012,edelsbrunner_computational_2010,kaczynski_computational_2004,otter_roadmap_2017}.

Topological data analysis is evolving rapidly. One particular direction of
evolution concerns a change in perspective from working with cell complexes as
topological spaces in and of themselves to focusing instead on data over a cell
complex --- viewing the cell complex as a base on which data to be investigated
resides. For example, one can consider scalar-valued data over cell complexes,
as occurs in weighted networks and complexes; or sensor data, as occurs in
target enumeration problems~\cite{curry_euler_2012}. Richer data involves
vector spaces and linear transformations, as with recent work in
cryo-EM~\cite{hadani_representation_2011} and synchronization
problems~\cite{bandeira_convex_2015}. Recent work in TDA points to a 
generalization of these and related data structures over
topological spaces. This is the theory of \style{sheaves}.

We will work exclusively with cellular sheaves~\cite{curry_sheaves_2014}. Fix a
(regular, locally finite) cell complex --- a triangulated surface will suffice
for purposes of imagination. A cellular sheaf of vector spaces is, in essence,
a data structure on this domain, assigning local data (in the form of vector
spaces) to cells and compatibility relations (linear transformations) between cells
of incident ascending dimension. These structure maps send data over vertices to data over
incident edges, data over edges to data over incident 2-cells, etc. As
a trivial example, the constant sheaf assigns a rank-one vector space to each
cell and identity isomorphisms according to boundary faces. More interesting is
the cellular analogue of a vector bundle: a cellular sheaf which assigns a
fixed vector space of dimension $n$ to each cell and isomorphisms as linear
transformations (with specializations to $O(n)$ or $SO(n)$ as desired).

The data assigned to a cellular sheaf naturally arranges into a cochain
complex graded by dimension of cells. 
As such, cellular sheaves possess a Laplacian that specializes to the
graph Laplacian and the Hodge Laplacian for the constant sheaf. For cellular
sheaves of real vector spaces, a spectral theory --- an examination of the
eigenvalues and eigenvectors of the sheaf Laplacian --- is natural, motivated,
and, to date, unexamined apart from a few special cases (see
\S\ref{sec:prior}).

This paper sketches an emerging spectral theory for cellular sheaves. Given the
motivation as a generalization of spectral graph theory, we will often
specialize to cellular sheaves over a 1-dimensional cell complex (that is, a
graph, allowing when necessary multiple edges between a pair of vertices). This is
mostly for the sake of simplicity and initial applications, as zero- and
one-dimensional homological invariants are the most readily applicable.
However, as the theory is general, we occasionally point to higher-dimensional
side-quests.

The plan of this paper is as follows. In \S\ref{sec:preliminaries}, we cover
the necessary topological and algebraic preliminaries, including definitions of
cellular sheaves. Next, \S\ref{sec:defns} gives definitions of the various
matrices involved in the extension of spectral theory to cellular sheaves.
Section~\ref{sec:harmonicity} uses these to explore issues related to harmonic
functions and cochains on sheaves. In \S\ref{sec:spectralbasics}, we extend
various elementary results from spectral graph theory to cellular sheaves. The
subsequent two sections  treat more sophisticated topics, effective resistance
(\S\ref{sec:resistance}) and the Cheeger inequality (\S\ref{sec:cheeger}), for
which we have some preliminary results. We conclude with outlines of potential
applications for the theory in \S\ref{sec:applications} and directions for
future inquiry in \S\ref{sec:further}.

The results and applications we sketch are at the beginnings of the subject,
and a great deal more in way of fundamental and applied work remains.

This paper has been written in order to be readable without particular
expertise in algebraic topology beyond the basic ideas of cellular homology and
cohomology. Category-theoretic terminology is used sparingly and for concision.
Given the well-earned reputation of sheaf theory as difficult for the
non-specialist, we have provided an introductory section with terminology and
core concepts, noting that much more is available in the
literature~\cite{bredon_sheaf_1997,kashiwara_sheaves_1990}. Our recourse to the
cellular theory greatly increases simplicity, readability, and applicability,
while resonating with the spirit of spectral graph theory. There are abundant
references available for the reader who requires more information on algebraic
topology~\cite{hatcher_algebraic_2001}, applications thereof~\cite{edelsbrunner_computational_2010,ghrist_elementary_2014}, and cellular
sheaf theory~\cite{curry_sheaves_2014,ghrist_elementary_2014}.

% ==========================================================================================
\section{Preliminaries}\label{sec:preliminaries}
% ==========================================================================================

% -------------------------------------------------------------------------------------
\subsection{Cell Complexes}\label{sec:cellcomp}
% -------------------------------------------------------------------------------------
\begin{definition}
A \style{regular cell complex} is a topological space $X$ with a partition into
subspaces $\{X_\alpha\}_{\alpha \in P_X}$ satisfying the following conditions:
\begin{enumerate}
\item
For each $x \in X$, every sufficiently small neighborhood of $x$ intersects finitely many $X_\alpha$.
\item
For all $\alpha,\beta$, $\overline{X_\alpha} \cap X_\beta \neq \varnothing$ only if $X_\beta \subseteq \overline{X_\alpha}$.
\item
Every $X_\alpha$ is homeomorphic to $\R^{n_\alpha}$ for some $n_\alpha$.
\item
For every $\alpha$, there is a homeomorphism of a closed ball in
$\R^{n_\alpha}$ to $\overline{X_\alpha}$ that maps the interior of the ball
homeomorphically onto $X_\alpha$.
\end{enumerate}
\end{definition}
Condition (2) implies that the set $P_X$ has a poset structure, given by $\beta
\leq \alpha$ iff $X_\beta \subseteq \overline{X_\alpha}$. This is known as the
\style{face poset} of $X$. The regularity condition (4) implies that all
topological information about $X$ is encoded in the poset structure of $P_X$.
For our purposes, we will identify a regular cell complex with its face poset,
writing the incidence relation $\beta\,\face\,\alpha$. The class of posets that
arise in this way can be characterized
combinatorially~\cite{bjorner_posets_1984}. For our purposes, a morphism of
cell complexes is a morphism of posets between their face incidence posets that
arises from a continuous map between their associated topological spaces. In
particular, morphisms of simplicial and cubical complexes will qualify as
morphisms of regular cell complexes.

The class of regular cell complexes includes simplicial complexes, cubical
complexes, and so-called multigraphs (as 1-dimensional cell complexes). As
nearly every space that can be characterized combinatorially can be represented
as a regular cell complex, these will serve well as a default class of spaces
over which to develop a combinatorial spectral theory of sheaves. We note that
the spectral theory of complexes has heretofore been largely restricted to the
study of simplicial complexes~\cite{schaub_random_2018}. A number of our
results will specialize to results about the spectra of Hodge Laplacians of regular cell complexes
by restricting to the constant sheaf.

A few notions associated to cell complexes will be useful. 
\begin{definition}
    The $k$-skeleton of a cell complex $X$, denoted
$X^{(k)}$,  is the subcomplex of $X$ consisting of cells of dimension at most $k$. 
\end{definition}
\begin{definition}
    Let $\sigma$ be a cell of a regular cell complex $X$. The star of $\sigma$,
    denoted $\st(\sigma)$, is the set of cells $\tau$ such that $\sigma \face
    \tau$. 
\end{definition}
Topologically, $\st(\sigma)$ is the smallest open collection of cells
containing $\sigma$, a role we might denote as the ``smallest cellular
neighborhood'' of $\sigma$. Stars serve an important purpose in giving combinatorial
analogues of topological notions for maps. For instance, a morphism $f:X \to Y$
of cell complexes may be locally injective as defined on the topological spaces.
Topologically, the condition for local injectivity is simply that every point in
$X$ have a neighborhood on which $f$ is injective. Translating this to cell
complexes, we require that for every cell $\sigma \in X$, $f$ is injective on
$\st(\sigma)$.

Topological continuity ensures that the preimage of a star $\st(\sigma)$ under a cell morphism
$f: X \to Y$ is a union of stars; if $f$ is locally injective, we see that it
must be a disjoint union of stars. A locally injective map is, further, a
covering map if on each component of $f^{-1}(\st(\sigma))$, $f$ is an
isomorphism. That is, the fiber of a star consists of a disjoint union of copies
of that star.
% -------------------------------------------------------------------------------------
\subsection{Cellular Sheaves}\label{sec:celsh}
% -------------------------------------------------------------------------------------
Let $X$ be a regular cell complex. A cellular sheaf attaches data spaces to the
cells of $X$ together with relations that specify when assignments to these data
spaces are consistent.

\begin{definition}
A \style{cellular sheaf} of vector spaces on a regular cell complex $X$ is an
assignment of a vector space $\Fc(\sigma)$ to each cell $\sigma$ of $X$
together with a linear transformation $\Fc_{\sigma\face\tau}\colon \Fc(\sigma)
\to \Fc(\tau)$ for each incident cell pair $\sigma\,\face\,\tau$. These must
satisfy both an identity relation $\Fc_{\sigma\face\sigma}=\id$ and the
composition condition:
\[
     \rho\,\face\,\sigma\,\face\,\tau
  ~~\Rightarrow~~
  \Fc_{\rho\face\tau} = \Fc_{\sigma\face\tau}\circ\Fc_{\rho\face\sigma}.
\]
The vector space $\Fc(\sigma)$ is called the \style{stalk} of $\Fc$ at
$\sigma$. The maps $\Fc_{\sigma\face\tau}$ are called the \style{restriction
maps}.
\end{definition}

For experts, this definition at first seems only reminiscent of the notion of
sheaves familiar to topologists. The depth of the relationship is explained in
detail in~\cite{curry_sheaves_2014}, but the essence is this: the data of a
cellular sheaf on $X$ specifies spaces of local sections on a cover of $X$
given by open stars of cells. This translates in two different ways into a
genuine sheaf on a topological space. One may either take the Alexandrov
topology on the face incidence poset of the complex, or one may view the open
stars of cells and their natural refinements a basis for the topology of $X$.
There then exists a natural completion of the data specified by the cellular
sheaf to a constructible sheaf on $X$.

One may compress the definition of a cellular sheaf to the following: If $X$ is
a regular cell complex with face incidence poset $P_X$, viewed as a category, a
cellular sheaf is a functor $\Fc\colon P_X \to \Vect_\Field$ to the category of
vector spaces over a field $\Field$.

\begin{definition}
Let $\Fc$ be a cellular sheaf on $X$. A \style{global section} $x$ of $\Fc$ is
a choice $x_\sigma \in \Fc(\sigma)$ for each cell $\sigma$ of $X$ such that
$x_\tau = \Fc_{\sigma\face\tau}x_\sigma$ for all $\sigma\,\face\,\tau$. The
space of global sections of $\Fc$ is denoted $\Gamma(X;\Fc)$.
\end{definition}

Perhaps the simplest sheaf on any complex is the constant sheaf with stalk
$\Vspace$, which we will denote $\constsh{\Vspace}$. This is the sheaf with all
stalks equal to $\Vspace$ and all restriction maps equal to the identity.

% -------------------------------------------------------------------------------------
\subsubsection{Cosheaves}\label{sec:cosh}
% -------------------------------------------------------------------------------------
In many situations it is more natural to consider a dual construction to a
cellular sheaf. A \style{cellular cosheaf} preserves stalk data but reverses
the direction of the face poset, and with it, the restriction maps.
\begin{definition}
A cellular cosheaf of vector spaces on a regular cell complex $X$ is an
assignment of a vector space $\Fc(\sigma)$ to each cell $\sigma$ of $X$
together with linear maps $\Fc_{\sigma\face\tau}\colon \Fc(\tau) \to
\Fc(\sigma)$ for each incident cell pair $\sigma\,\face\,\tau$ which satisfies
the identity ($\Fc_{\sigma\face\sigma}=\id$) and composition condition:
\[
     \rho\,\face\,\sigma\,\face\,\tau
  ~~\Rightarrow~~
  \Fc_{\rho\face\tau} = \Fc_{\rho\face\sigma}\circ\Fc_{\sigma\face\tau} .
\]
\end{definition}
More concisely, a cellular cosheaf is a functor $P_X^{\text{op}} \to
\Vect_\Field$. The contravariant functor $\Hom(\bullet,\Field):
\Vect_\Field^\text{op} \to \Vect_\Field$ gives every cellular sheaf $\Fc$ a
dual cosheaf $\hat{\Fc}$ whose stalks are $\Hom(\Fc(\sigma),\Field)$.

% -------------------------------------------------------------------------------------
\subsubsection{Homology and Cohomology}\label{sec:homcohom}
% -------------------------------------------------------------------------------------
The cells of a regular cell complex have a natural grading by dimension. By
regularity of the cell complex, this grading can be extracted from the face
incidence poset as the height of a cell in the poset. This means that a
cellular sheaf has a graded vector space of cochains
\[C^k(X;\Fc) = \bigoplus_{\dim(\sigma) = k} \Fc(\sigma).\]

To develop this into a chain complex, we need a boundary operator and a notion
of orientation --- a signed incidence relation on $P_X$. This is a map
$[\bullet:\bullet]: P_X \times P_X \to \{0,\pm 1\}$ satisfying the following
conditions:
\begin{enumerate}
\item
If $[\sigma:\tau] \neq 0$, then $\sigma\face\tau$ and there are no cells between $\sigma$ and $\tau$ in the incidence poset.
\item For any $\sigma\face\tau$, $\sum_{\gamma \in P_X} [\sigma:\gamma][\gamma:\tau] = 0$.
\end{enumerate}
Given a signed incidence relation on $P_X$, there exist coboundary maps $\delta^k: C^k(X;\Fc) \to C^{k+1}(X;\Fc)$. These are given by the formula
\[\delta^k|_{\Fc(\sigma)} = \sum_{\dim(\tau) = k+1} [\sigma:\tau]\Fc_{\sigma\face\tau},\]
or equivalently,
\[(\delta^kx)_\tau = \sum_{\dim(\sigma) = k} [\sigma:\tau] \Fc_{\sigma\face\tau} (x_\sigma).\]
Here we use subscripts to denote the value of a cochain in a particular stalk;
that is, $x_\sigma$ is the value of the cochain $x$ in the stalk $\Fc(\sigma)$.

It is a simple consequence of the properties of the incidence relation and the
commutativity of the restriction maps that $\delta^k \circ \delta^{k-1} = 0$, so these
coboundary maps define a cochain complex and hence a cohomology theory for
cellular sheaves. In particular, $H^0(X;\Fc)$ is naturally isomorphic to
$\Gamma(X;\Fc)$, the space of global sections. An analogous construction
defines a homology theory for cosheaves. Cosheaf homology may be thought of as
dual to sheaf cohomology in a Poincar\'{e}-like sense. That is, frequently the
natural analogue of degree zero sheaf cohomology is degree $n$ cosheaf
homology. A deeper formal version of this fact, exploiting an equivalence of
derived categories, may be found in~\cite[ch. 12]{curry_sheaves_2014}.

There is a relative version of cellular sheaf cohomology. Let $A$ be a
subcomplex of $X$. There is a natural subspace of $C^k(X;\Fc)$ consisting of
cochains which vanish on stalks over cells in $A$. The coboundary of a cochain
which vanishes on $A$ also vanishes on $A$, since any cell in $A^{(k+1)}$ has
only cells in $A^{(k)}$ on its boundary. We therefore get a subcomplex
$C^\bullet(X,A;\Fc)$ of $C^\bullet(X;\Fc)$. The cohomology of this subcomplex
is the relative sheaf cohomology $H^\bullet(X,A;\Fc)$. The natural maps between
these spaces of cochains constitute a short exact sequence of complexes
\[
    0 \to C^\bullet(X,A;\Fc) \to C^\bullet(X;\Fc) \to C^\bullet(A;\Fc) \to 0,
\]
from which a long exact sequence for relative sheaf cohomology arises:
\[
    0 \to H^0(X,A;\Fc) \to H^0(X;\Fc) \to H^0(A;\Fc) \to H^1(X,A;\Fc) \to \cdots
\]
% -------------------------------------------------------------------------------------
\subsubsection{Sheaf Morphisms}
\label{sec:morph}
% -------------------------------------------------------------------------------------
\begin{definition}
If $\Fc$ and $\Gc$ are sheaves on a cell complex $X$, a sheaf \style{morphism}
$\varphi: \Fc \to \Gc$ is a collection of maps $\varphi_\sigma: \Fc(\sigma) \to
\Gc(\sigma)$ for each cell $\sigma$ of $X$, such that for any
$\sigma\face\tau$, $\varphi_\tau \circ \Fc_{\sigma\face\tau} =
\Gc_{\sigma\face\tau} \circ \varphi_\sigma$. Equivalently, all diagrams of the
following form commute:
\[
\begin{tikzcd}
  \Fc(\sigma) \arrow[d,"\Fc_{\sigma \face\tau}"'] \arrow[r,"\varphi_\sigma"] & \Gc(\sigma) \arrow[d,"\Gc_{\sigma\face \tau}"] \\
  \Fc(\tau) \arrow[r,"\varphi_\tau"] & \Gc(\tau)
\end{tikzcd}
\]
\end{definition}
This commutativity condition assures that a sheaf morphism $\varphi: \Fc \to
\Gc$ induces maps $\varphi^k: C^k(X;\Fc) \to C^k(X;\Gc)$ which commute with the
coboundary maps, resulting in the induced maps on cohomology  $H^k \varphi:
H^k(X;\Fc) \to H^k(X;\Gc)$.

% -------------------------------------------------------------------------------------
\subsubsection{Sheaf Operations}
\label{sec:ops}
% -------------------------------------------------------------------------------------
There are several standard operations that act on sheaves to produce new sheaves.
\begin{definition}[Direct sum]
If $\Fc$ and $\Gc$ are sheaves on $X$, their \style{direct sum} $\Fc \oplus
\Gc$ is a sheaf on $X$ with $(\Fc \oplus \Gc)(\sigma) = \Fc(\sigma) \oplus
\Gc(\sigma)$. The restriction maps are $(\Fc\oplus\Gc)_{\sigma\face\tau} =
\Fc_{\sigma\face \tau} \oplus \Gc_{\sigma\face\tau}$.
\end{definition}

\begin{definition}[Tensor product]
If $\Fc$ and $\Gc$ are sheaves on $X$, their \style{tensor product} $\Fc
\otimes \Gc$ is a sheaf on $X$ with $(\Fc \otimes \Gc)(\sigma) = \Fc(\sigma)
\otimes \Gc(\sigma)$. The restriction maps are $(\Fc \otimes
\Gc)_{\sigma\face\tau} = \Fc_{\sigma\face \tau} \otimes \Gc_{\sigma\face\tau}$.
\end{definition}

\begin{definition}[Pullback]
If $f:X \to Y$ is a morphism of cell complexes and $\Fc$ is a sheaf on $Y$, the
\style{pullback} $f^*\Fc$ is a sheaf on $X$ with $f^*\Fc(\sigma) =
\Fc(f(\sigma))$ and $(f^*\Fc)_{\sigma\face \tau} = \Fc_{f(\sigma)\face
f(\tau)}$.
\end{definition}

\begin{definition}[Pushforward]
The full definition of the pushforward of a cellular sheaf is somewhat more
categorically involved than the previous constructions. If $f: X \to Y$ is a
morphism of cell complexes and $\Fc$ is a sheaf on $X$, the \style{pushforward}
$f_*\Fc$ is a sheaf on $Y$ with stalks $f_*\Fc(\sigma)$ given as the limit
$\lim_{\sigma \face f(\tau)} \Fc(\tau)$. The restriction maps are induced by
the restriction maps of $\Fc$, since whenever $\sigma \face \sigma'$, the cone
for the limit defining $f_*\Fc(\sigma)$ contains the cone for the limit
defining $f_*\Fc(\sigma')$, inducing a unique map $f_*\Fc(\sigma) \to
f_*\Fc(\sigma')$.

In this paper, we will mainly work with pushforwards over locally injective cell
maps, that is, those whose geometric realizations are locally injective (see
\S\ref{sec:cellcomp}). If $f: X \to Y$
is locally injective, every cell $\sigma \in X$ maps to a cell of the
same dimension, and for every cell $\sigma \in Y$, $f^{-1}(\st(\sigma))$ is a
disjoint union of subcomplexes, each of which maps injectively to $Y$. In this
case, $f^*\Fc(\sigma) \simeq \bigoplus_{\sigma' \in f^{-1}(\sigma)}
\Fc(\sigma')$, and $(f^*\Fc)_{\sigma\face\tau} = \bigoplus_{(\sigma' \face
\tau') \in f^{-1}(\sigma\face\tau)} \Fc_{\sigma'\face\tau'}$.
This computational formula in fact holds more generally, if the stars of cells
in $f^{-1}(\sigma)$ are disjoint.

\end{definition}

Those familiar with the definitions of pushforward and pullback for sheaves
over topological spaces will note a reversal of fates when we define sheaves
over cell complexes. Here the pullback is simple to define, while the
pushforward is more involved. This complication arises because cellular sheaves
are in a sense defined pointwise rather than over open sets.

% ==========================================================================================
\section{Definitions}\label{sec:defns}
% ==========================================================================================

% -------------------------------------------------------------------------------------
\subsection{Weighted Cellular Sheaves}\label{sec:weighted}
% -------------------------------------------------------------------------------------
Let $\Field = \R$ or $\C$. A weighted cellular sheaf is a cellular sheaf with
values in $\Field$-vector spaces where the stalks have additionally been given
an inner product structure. Adding the condition of completeness to the stalks,
one may view this as a functor $P_X \to \Hilb_\Field$, where $\Hilb_\Field$ 
is the category whose objects are Hilbert
spaces over $\Field$ and whose morphisms are (bounded) linear maps.

The inner products on stalks of $\Fc$ extend by the orthogonal direct sum to
inner products on $C^k(X;\Fc)$, making these Hilbert spaces as well.  The
canonical inner products on direct sums and subspaces of Hilbert spaces give the
direct sum and tensor product of weighted cellular sheaves weighted structures.
Similarly, the pullbacks and pushforwards (over locally injective maps) of a
weighted sheaf have canonical weighted structures given by their computational
formulae in \S\ref{sec:ops}.

Every morphism $T:V \to W$ between Hilbert spaces admits an adjoint map $T^*:W
\to V$, determined by the property that for all $v \in V, w \in W$, $\ip{w,Tv}
= \ip{T^*w,v}$. One may readily check that $(T^*)^* = T$. This fact gives the
category $\Hilb_\Field$ a \emph{dagger structure}, that is, a contravariant
endofunctor $\dagger$ (here the adjoint operation ${}^*$) which acts as the identity on objects and squares to the
identity. In a dagger category, the notion of unitary isomorphisms makes sense:
they are the invertible morphisms $T$ such that $T^\dagger = T^{-1}$. 

The dagger structure of $\Hilb_\Field$ introduces some categorical subtleties
into the study of weighted cellular sheaves. The space of global sections of a
cellular sheaf is defined in categorical terms as the limit of the functor $X
\to \Vect$ defining the sheaf. This defines the space of global sections up to
unique isomorphism. We might want a weighted space of global sections to be a
sort of limit in $\Hilb_\Field$ which is defined up to unique \emph{unitary}
isomorphism. This is the notion of \emph{dagger limit}, recently studied in~\cite{heunen_limits_2019}. 
Unfortunately, this work showed that $\Hilb_\Field$
does not have all dagger limits; in particular, pullbacks over spans of
noninjective maps do not exist. As a result, there is
no single canonical way to define an inner product on the space of global
sections of a cellular sheaf $\Fc$. There are two approaches that seem most
natural, however. One is to view the space of global sections of $\Fc$ as $\ker
\delta^0_\Fc$ with the natural inner product given by inclusion into
$C^0(X;\Fc)$. The other is to view global sections as lying in
$\bigoplus_{\sigma} \Fc(\sigma)$. We will generally take the view that global
sections are a subspace of $C^0(X;\Fc)$; that is, we will weight
$\Gamma(X;\Fc)$ by its canonical isomorphism with $\mathcal{H}^0(X;\Fc)$, as
defined in \S\ref{sec:sheaflap}.

The dagger structure on $\Hilb_\Field$ gives a slightly different way to
construct a dual cosheaf from a weighted cellular sheaf $\Fc$. Taking the
adjoint of each restriction map reverses their directions and hence yields a
cosheaf with the same stalks as the original sheaf. From a categorical
perspective, this amounts to composing the functor $\Fc$ with the dagger
endofunctor on $\Hilb_\Field$. When stalks are finite
dimensional, this dual cosheaf is isomorphic to the cosheaf $\hat{\Fc}$ defined
in \S\ref{sec:cosh} via the dual vector spaces of stalks.  In this situation,
we have an isomorphism between the stalks of $\Fc$ and its dual cosheaf.  This
is reminiscent of the bisheaves recently introduced by MacPherson and
Patel~\cite{macpherson_persistent_2018}. However, the structure maps
$\Fc(\sigma) \to \hat{\Fc}(\sigma)$ will rarely commute with the restriction
and extension maps as required by the definition of the bisheaf --- this only
holds in general if all restriction maps are unitary. The bisheaf construction
is meant to give a generalization of local systems, and as such fits better
with our discussion of discrete vector bundles in \S\ref{sec:discvect}.

% -------------------------------------------------------------------------------------
\subsection{The Sheaf Laplacian}\label{sec:sheaflap}
% -------------------------------------------------------------------------------------
Given a chain complex of Hilbert spaces $C^0 \to C^1 \to
\cdots$ we can construct the Hodge Laplacian $\Delta = (\delta + \delta^*)^2 =
\delta^*\delta + \delta \delta^*$. This operator is naturally graded into
components $\Delta^k: C^k \to C^k$, with $\Delta^k = (\delta^k)^*\delta^k +
\delta^{k-1}(\delta^{k-1})^*$. This operator can be further separated into
\style{up-} (coboundary) and \style{down-} (boundary) Laplacians $\Delta_+^k =
(\delta^k)^*\delta^k$ and $\Delta_-^k = \delta^{k-1}(\delta^{k-1})^*$
respectively.

A key observation is that on a
finite-dimensional Hilbert space, $\ker \delta^* = (\im \delta)^\perp$. For if
$\delta^* x = 0$, then for all $y$, $0 = \ip{\delta^* x, y} = \ip{x,\delta y}$,
so that $x \perp \im \delta$. This allows us to express the kernels and images
necessary to compute cohomology purely in terms of kernels. This is the content
of the central theorem of discrete Hodge theory:
\begin{theorem}
    Let $C^0 \to C^1 \to \cdots$ be a chain complex of finite-dimensional Hilbert spaces, with Hodge Laplacians $\Delta^k$. Then $\ker \Delta^k \cong H^k(C^\bullet)$.
\end{theorem}
\begin{proof}
    By definition, $H^k(C^\bullet) = \ker \delta^k /\img \delta^{k-1}$. In a
    finite dimensional Hilbert space, $\ker \delta^k / \img \delta^{k-1}$ is
    isomorphic to the orthogonal complement of $\img \delta^{k-1}$ in $\ker
    \delta^k$, which we may write $(\ker \delta^k) \cap (\img
    \delta^{k-1})^\perp = (\ker \delta^k) \cap (\ker (\delta^{k-1})^*)$.  So it
    suffices to show that $\ker \Delta^k = (\ker \delta^k) \cap (\ker
    (\delta^{k-1})^*)$. Note that $\ker \delta^k = \ker (\delta^k)^*\delta^k =
    \ker \Delta^k_+$ and similarly for $\Delta^k_-$. So we need to show that
    $\ker (\Delta^k_+ + \Delta^k_-) = \ker \Delta^k_+ \cap \ker \Delta^k_-$,
    which will be true if $\im \Delta^k_+ \cap \im \Delta^k_- = 0$. But this is
    true because $\im \Delta^k_+ = \im (\delta^k)^* = (\ker\delta^k)^\perp$ and
    $\im \Delta^k_- = \im \delta^{k-1} \subseteq \ker \delta^k$.
\end{proof}

The upshot of this theorem is that the kernel of $\Delta^k$ gives a set of
canonical representatives for elements of $H^k(C^\bullet)$. This is commonly
known as the space of \style{harmonic cochains}, denoted $\mathcal
H^k(C^\bullet)$. In particular, the proof above implies that there is an
orthogonal decomposition $C^k = \mathcal H^k \oplus \im \delta^{k-1} \oplus \im
(\delta^{k})^*$.

When the chain complex in question is the complex of cochains for a weighted
cellular sheaf $\Fc$, the Hodge construction produces the sheaf Laplacians. The
Laplacian which is easiest to study and most immediately interesting is the
degree-0 Laplacian, which is a generalization of the graph Laplacian. We can
represent it as a symmetric block matrix with blocks indexed by the vertices of
the complex. The entries on the diagonal are $\Delta^0_{v,v} = \sum_{v\face e}
\Fc_{v\face e}^*\Fc_{v\face e}$ and the entries on the off-diagonal are
$\Delta^0_{u,v} = -\Fc_{u\face e}^*\Fc_{v\face e}$, where $e$ is the edge
between $v$ and $u$.  Laplacians of other degrees have similar block
structures.

The majority of results in combinatorial spectral theory have to do with
up-Laplacians. We will frequently denote these $L^k$ by analogy with spectral
graph theory, where $L$ typically denotes the (non-normalized) graph Laplacian.
In particular, we will further elide the index $k$ when $k=0$, denoting the
graph sheaf Laplacian by simply $L$. A subscript will be added when necessary
to identify the sheaf, \eg~$L_\Fc$ or $\Delta_\Fc^k$.

\begin{figure}[hb]
    {
    \centering
\includegraphics[width=4in]{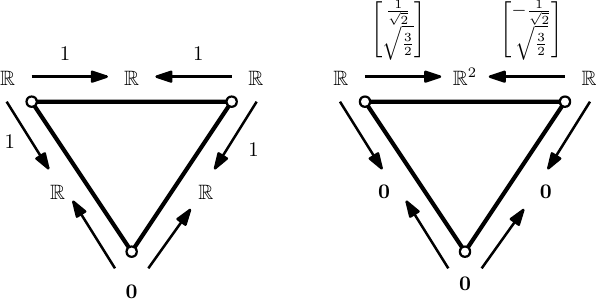}
\caption{Two nonisomorphic sheaves with the same Laplacian.}\label{fig:nonisomorphic}
}
\end{figure}

Weighted labeled graphs are in one-to-one correspondence with graph Laplacians.
The analogous statement is not true of sheaves on a graph. For instance, the
sheaves in Figure~\ref{fig:nonisomorphic} have coboundary maps with matrix
representations
\[
    \begin{bmatrix}
        1 & -1 \\
        1 & 0 \\
        0 & 1
    \end{bmatrix}
    \quad\text{and}\quad
    \begin{bmatrix}
        \frac{1}{\sqrt{2}} & \frac{1}{\sqrt{2}} \\
        \sqrt{\frac{3}{2}} & -\sqrt{\frac{3}{2}}
    \end{bmatrix},
\]
which means that the Laplacian for each is equal to
\[
    \begin{bmatrix}
        2 & -1 \\ 
        -1 & 2
    \end{bmatrix}.
\]
However, these sheaves are not unitarily isomorphic, as can be seen immediately
by checking the stalk dimensions. More pithily, one cannot hear the shape of a
sheaf. One source of the lossiness in the sheaf Laplacian representation 
is that restriction maps may be the zero morphism,
effectively allowing for edges that are only attached to one vertex. More
generally, restriction maps may fail to be full rank, which means that it is
impossible to identify the dimensions of edge stalks from the Laplacian.

% -------------------------------------------------------------------------------------
\subsubsection{Harmonic Cochains}\label{sec:harmonics}
% -------------------------------------------------------------------------------------
The elements of $\ker \Delta^k = \mathcal H^k$ are known as harmonic
$k$-cochains. More generally, a $k$-cochain may be harmonic on a subcomplex:

\begin{definition}
A $k$-cochain $x$ of a sheaf $\Fc$ on a cell complex $X$ is \style{harmonic} on
a set $S$ of $k$-cells if $(\Delta^k_\Fc x)|_S = 0$.
\end{definition}
When $k = 0$ and $\Fc$ is the constant sheaf (\ie, in spectral graph theory),
this can be expressed as a local averaging property: For each $v \in S$, $x_v =
\frac{1}{d_v} \sum_{u \sim v} x_u$, where $\sim$ indicates adjacency and $d_v$
is the degree of the vertex $v$.

% -------------------------------------------------------------------------------------
\subsubsection{Identifying Sheaf Laplacians}\label{sec:ident}
% -------------------------------------------------------------------------------------
Given a regular cell complex $X$ and a choice of dimension for each stalk, one
can identify the collection of matrices which arise as coboundary maps for a
sheaf on $X$ as those matrices satisfying a particular block sparsity pattern.
This sparsity pattern controls the number of nonzero blocks in each row of the
matrix. Restricting to $\delta^0$, we get a matrix whose rows have at most two
nonzero blocks. The space of matrices which arise as sheaf Laplacians is then
the space of matrices which have a factorization $L = \delta^*\delta$, where
$\delta$ is a matrix satisfying this block sparsity condition. Boman et al.\
studied this class of matrices when the blocks have size $1\times 1$,
calling them matrices of \style{factor width} two~\cite{boman_factor_2005}.
They showed that this class coincides with the class of symmetric generalized
diagonally dominant matrices, those matrices $L$ for which there exists a
positive diagonal matrix $D$ such that $DLD$ is diagonally dominant. 
Indeed, the fact that sheaves on graphs are not in general determined by their
Laplacians is in part a consequence of the nonuniqueness of width-two
factorizations.

% -------------------------------------------------------------------------------------
\subsection{Approaching Infinite-Dimensional Laplacians}\label{sec:infdim}
% -------------------------------------------------------------------------------------

The definitions given in this paper are adapted to the case of sheaves of
finite dimensional Hilbert spaces over finite cell complexes. Relaxing these
finiteness constraints introduces new complications. 

The spaces of cochains
naturally acquire inner products by taking the Hilbert space direct sum. These
are not the same
as taking the algebraic direct sum or product of stalks. However, there is a
sequence of inclusions of complexes
\[
    C^\bullet_c(X;\Fc) \subseteq L^2C^\bullet(X;\Fc) \subseteq C^\bullet(X;\Fc)
\]
inducing algebraic maps between the corresponding compactly supported, $L^2$, and standard sheaf cohomology theories.

The theory of abstract complexes of possibly infinite-dimensional Hilbert spaces has
been developed in~\cite{bruning_hilbert_1992}. This paper explains conditions
for the spaces of harmonic cochains of a complex to be isomorphic with the
algebraic cohomology of the complex. A particularly nice condition is that the
complex have finitely generated cohomology, which implies that the total
coboundary map is a Fredholm operator. More generally, if the images of the
coboundary and its adjoint are closed, the spaces of harmonic cochains will be
isomorphic to the cohomology.

Further issues arise when we consider the coboundary maps $\delta^k$. 
For spectral purposes, it is in general desirable for these to be 
bounded linear maps, for which we must make some further
stipulations. Sufficient conditions for coboundary maps to be bounded are as follows:
\begin{proposition}
    Let $\Fc$ be a sheaf of Hilbert spaces on a cell complex $X$. Suppose that
    there exists some $M_k$ such that for every pair of cells $\sigma \face \tau$
    with $\dim \sigma = k$ and $\dim \tau = k+1$, $\norm{\Fc_{\sigma \face
    \tau}} \leq M_k$. Further suppose that every $k$-cell in $X$ has at most
    $d^k$ cofaces of dimension $k+1$, and every $(k+1)$-cell in $X$ has at most
    $d_{k+1}$ faces of dimension $k$. Then $\delta^k_\Fc$ is a bounded linear
    operator.
\end{proposition}
\begin{proof}
    We compute:
    \begin{multline*}
    \norm{\delta^kx}^2 = \sum_{\dim \tau = k+1} \norm{(\delta^kx)_\tau}^2 \leq
    \sum_{\dim \tau = k+1} \left(\sum_{\sigma \face \tau} \norm{\Fc_{\sigma \face \tau}x_\sigma}\right)^2 \\
    \leq \sum_{\dim \tau = k+1} \left(\sum_{\sigma \face \tau} M_k \norm{x_\sigma}\right)^2 
    \leq M_k^2 \sum_{\dim \tau = k+1} d_{k+1}\sum_{\sigma \face \tau}
    \norm{x_\sigma}^2 \\
    = M_k^2 d_{k+1}\sum_{\dim \sigma = k} \sum_{\sigma \face \tau}
    \norm{x_\sigma}^2
    \leq M_k^2 d_{k+1}d^k \sum_{\dim \sigma = k} \norm{x_\sigma}^2
    = M_k^2 d_{k+1}d^k \norm{x}^2.
\end{multline*}
\end{proof}
If $\delta^k$ is bounded, its associated Laplacians $\Delta^k_+ =
(\delta^k)^*\delta^k$ and $\Delta^{k+1}_- = \delta^k(\delta^k)^*$ are also
bounded. As bounded self-adjoint operators, their spectral theory is relatively
unproblematic. Their spectra consist entirely of approximate eigenvalues, those
$\lambda$ for which there exists a sequence of unit vectors $\{x_k\}$ such that $\norm{\Delta^k_+x_k
- \lambda x_k} \to 0$. 

If $\delta^k$ is not just bounded, but compact, the Laplacian spectral theory
becomes even nicer. In this situation, the spectrum of $\Delta^k_+$ has no
continuous part, and hence consists purely of eigenvalues. An
appropriate decay condition on norms of restriction maps ensures compactness.
\begin{proposition}
    Let $\Fc$ be a sheaf of Hilbert spaces on a cell complex $X$. Suppose that
    for all $\sigma \face \tau$ with $\dim \sigma = k$ and
    $\dim \tau = k+1$, the restriction map $\Fc_{\sigma \face \tau}$ for is compact, 
    and further that $\sum_{\sigma \face \tau}
    \norm{\Fc_{\sigma \face \tau}} < \infty$. Then $\delta^k_\Fc$ is a compact
    linear operator.
\end{proposition}
\begin{proof}
    It is clear that $\delta^k$ cannot be compact if any one of its component restriction
    maps fails to be compact. Suppose first that all restriction maps are finite
    rank, and fix an ordering of $(k+1)$-cells of $X$, defining the orthogonal
    projection operators $P^i: C^{k+1}(X;\Fc) \to C^{k+1}(X;\Fc)$ sending stalks
    over $(k+1)$ cells of index greater than $i$ to zero. Then $P^i \delta^k$ is
    a finite-rank operator and
    \[ \norm{P^i \delta^k - \delta^k} \leq \sum_{j > i} \sum_{\sigma \face
    \tau_j}\norm{\Fc_{\sigma \face \tau_j}},\]
    which goes to zero as $i \to \infty$. In the case that the restriction maps
    are compact, pick an approximating sequence for each by finite rank maps and
    combine the two approximations.
\end{proof}

An important note is that when $C^k(X;\Fc)$ is infinite dimensional and
$\delta^k$ is compact with finite dimensional kernel, the eigenvalues of $\Delta^k_+$ will accumulate at zero.
This means that there will be no smallest nontrivial eigenvalue for such
Laplacians.

Most of the difficulties considered here already arise in the study of spectra of infinite graphs.
The standard Laplacian associated to an infinite graph is bounded but not
compact, while a proper choice of weights decaying at infinity makes it
compact.

The study of sheaves of arbitrary
Hilbert spaces on not-necessarily-finite cell complexes is interesting, and
indeed suggests itself in certain applications. However, for
the initial development and exposition of the theory, we have elected to focus on the (still
quite interesting) finite-dimensional case. This is sufficient for most
applications we have envisioned, and avoids the need for repeated
qualifications and restrictions. 

For the balance of this paper, we will assume that all cell complexes are finite
and all vector spaces are finite dimensional,
giving where possible proofs that generalize in some way to the
infinite-dimensional setting. Most results that do not explicitly require a
finite complex will extend quite directly to the case of sheaves with compact
coboundary operators. Proofs not relying on the Courant-Fischer theorem will
typically apply even to situations where coboundary operators are merely
bounded, although their conclusions may be somewhat weakened.  

% -------------------------------------------------------------------------------------
\subsection{The Normalized Laplacian and Weights}
\label{sec:norm}
% -------------------------------------------------------------------------------------
Many results in spectral graph theory rely on a normalized version of the
standard graph Laplacian, which is typically defined in terms of a rescaling of
the standard Laplacian. Let $D$ be the diagonal matrix whose nonzero entries
are the degrees of vertices; then the normalized Laplacian is $\mathcal L =
D^{-1/2}L D^{-1/2}$. This definition preserves the Laplacian as a symmetric
matrix, but it obscures the true meaning of the normalization. The normalized
Laplacian is the standard Laplacian with a different choice of weights for the
vertices. The matrix $D^{-1/2}L D^{-1/2}$ is similar to $D^{-1}L$, which is
self adjoint with respect to the inner product $\ip{x,y} = x^TDy$. In this
interpretation, each vertex is weighted proportionally to its degree. Viewing the
normalization process as a reweighting of cells leads to the natural definition
of normalized Laplacians for simplicial complexes given by Horak and
Jost~\cite{horak_spectra_2013}.

Indeed, following Horak and Jost's definition for simplicial complexes, we
propose the following extension to sheaves.
\begin{definition}\label{defn_normalized_sheaf}
Let $\Fc$ be a weighted cellular sheaf defined on a regular cell complex $X$.
We say $\Fc$ is \style{normalized} if for every cell $\sigma$ of $X$ and every
$x,y \in \Fc(\sigma) \cap (\ker \delta)^\perp$, $\ip{\delta x, \delta y} =
\ip{x,y}$.
\end{definition}

\begin{lemma} Given a weighted sheaf $\Fc$ on a finite-dimensional cell complex
$X$, it is always possible to reweight $\Fc$ to a normalized version.
\end{lemma}
\begin{proof} Note that if $X$ has dimension $k$, the normalization condition is
    trivially satisfied for all cells $\sigma$ of dimension $k$. Thus, starting
    at cells of dimension $k-1$, we recursively redefine the inner products on
    stalks. If $\sigma$ is a cell of dimension $k-1$, let $\Pi_\sigma$ be the
    orthogonal projection $\Fc(\sigma) \to \Fc(\sigma)\cap \ker \delta$. Then
    define the normalized inner product $\ip{\bullet,\bullet}^N_\sigma$ on
    $\Fc(\sigma)$ to be given by $\ip{x,y}^N_{\sigma} = \ip{\delta (\id -
        \Pi_\sigma)x,\delta (\id - \Pi_\sigma)y} + \ip{\Pi_\sigma x, \Pi_\sigma
    y}$. It is clear that this reweighted sheaf satisfies the condition of
    Definition~\ref{defn_normalized_sheaf} for cells of dimension $k$ and $k-1$.
    We may then perform this operation on cells of progressively lower dimension
    to obtain a fully normalized sheaf.  \end{proof}

Note that there is an important change of perspective here: we do not normalize
the Laplacian of a sheaf, but instead normalize the sheaf itself, or more
specifically, the inner products associated with each stalk of the sheaf.

If we apply this process to a sheaf $\Fc$ on a graph $G$, there is an immediate
interpretation in terms of the original sheaf Laplacian. Let $D$ be the block
diagonal of the standard degree 0 sheaf Laplacian, and note that for $x \perp \ker L$,
$\ip{x,Dx}$ is the reweighted inner product on $C^0(G;\Fc)$. In particular, the
adjoint of $\delta$ with respect to this inner product has the form $D^\dagger
\delta^T$, where $D^\dagger$ is the Moore-Penrose pseudoinverse of $D$, so that
the matrix form of the reweighted Laplacian with respect to this inner product
is $D^\dagger L$. Changing to the standard basis then gives $\mathcal L = D^{\dagger/2} L D^{\dagger/2}$.

% -------------------------------------------------------------------------------------
\subsection{Discrete Vector Bundles}\label{sec:discvect}
% -------------------------------------------------------------------------------------
A subclass of sheaves of particular interest are those where all restriction
maps are invertible.These sheaves have been the subject of significantly more
study than the general case, since they extend to locally constant sheaves on
the geometric realization of the cell complex. The Riemann-Hilbert
correspondence describes an equivalence between locally constant sheaves (or
cosheaves) on $X$, local systems on $X$, vector bundles on $X$ with a flat
connection, and representations of the fundamental groupoid of $X$. (See,
\eg,~\cite[ch. 5]{davis_lecture_2001} or~\cite{bass_local_2009} for a
discussion of some aspects of this correspondence.) When we represent a local
system by a cellular sheaf or cosheaf, we will call it a discrete vector
bundle.

One way to understand the space of $0$-cochains of a discrete vector bundle is
as representing a subspace of the sections of a geometric realization of the
associated flat vector bundle, defined by linear interpolation over
higher-dimensional cells. The coboundary map can be seen as a sort of
discretization of the connection, whose flatness is manifest in the fact that
$\delta^2 = 0$.

Discrete vector bundles have some subtleties when we study their Laplacians.
The sheaf-cosheaf duality corresponding to a local system, given by taking
inverses of restriction maps, is not in general the same as the duality induced
by an inner product on stalks. Indeed, these duals are only the same when the
restriction maps are unitary --- their adjoints must be their inverses.

The inner product on stalks of a cellular sheaf has two roles: it gives a
relative weight to vectors in each stalk, but via the restriction maps also
gives a relative weight to cells in the complex. This second role complicates
our interpretation of certain sorts of vector bundles. For instance, one might
wish to define an $O(n)$ discrete vector bundle on a graph to be a cellular
sheaf of real vector spaces where all
restriction maps are orthogonal. However, from the perspective of the degree-0
Laplacian, a uniform scaling of the inner product on an edge does not change
the orthogonality of the bundle, but instead in some sense changes the length
of the edge, or perhaps the degree of emphasis we give to discrepancies over
that edge. So a discrete $O(n)$-bundle should be one where the restriction maps
on each cell are {\em scalar multiples} of orthonormal maps.

That is, for each cell $\sigma$, we have a positive scalar $\alpha_\sigma$,
such that for every $\sigma\face\tau$, the restriction map
$\Fc_{\sigma\face\tau}$ is an orthonormal map times
$\alpha_{\tau}/\alpha_\sigma$. One way to think of this is as a scaling of the
inner product on each stalk of $\Fc$. Frequently, especially when dealing with
graphs, we set $\alpha_\sigma = 1$ when $\dim(\sigma) = 0$, but this is not
necessary. (Indeed, when dealing with the normalized Laplacian of a graph, we
have $\alpha_v = \sqrt{d_v}$.) 

The rationale for this particular definition is that in the absence of a basis,
inner products are not absolutely defined, but only in relation to maps in or
out of a space. Scaling the inner product on a vector space is meaningless
except in relation to a given collection of maps, which it transforms in a
uniform way.

As a special case of this definition, it will be useful to think about weighted
versions of the constant sheaf. These are isomorphic to the `true' constant
sheaf, but not unitarily so. Weighted constant sheaves on a graph are analogous
to weighted graphs. The distinction between the true constant sheaf and
weighted versions arises because it is often convenient to think of the
sections of a cellular sheaf as a subspace of $C^0(X;\Fc)$. As a result, we
often only want our sections to be constant on $0$-cells, allowing for
variation up to a scalar multiple on higher-dimensional cells. This notion will
be necessary in \S\ref{sec:approx_sheaf} when we discuss approximations of
cellular sheaves.

% -------------------------------------------------------------------------------------
\subsection{Comparison with Previous Constructions}\label{sec:prior}
% -------------------------------------------------------------------------------------
Friedman, in~\cite{friedman_sheaves_2015}, gave a definition of a
sheaf\footnote{In our terminology, Friedman's sheaves are cellular cosheaves.}
on a graph, developed a homology theory, and suggested constructing sheaf
Laplacians and adjacency matrices.  The suggestion that one might develop a
spectral theory of sheaves on graphs has remained until now merely a
suggestion.

The graph connection Laplacian, introduced by Singer and Wu
in~\cite{singer_vector_2012}, is simply the sheaf Laplacian of an $O(n)$-vector
bundle over a graph. This construction has attracted significant interest from
a spectral graph theory perspective, including the development of a
Cheeger-type inequality~\cite{bandeira_cheeger_2013} and a study of random
walks and sparsification~\cite{chung_ranking_2012}. Connection Laplacian
methods have proven enlightening in the study of synchronization problems.
Others have approached the study of vector bundles, and in particular line
bundles, over graphs without reference to the connection Laplacian, studying
analogues of spanning trees and the Kirchhoff
theorems~\cite{kenyon_spanning_2011,catanzaro_kirchhoffs_2013}. Other work on
discrete approximations to connection Laplacians of manifolds has analyzed
similar matrices~\cite{mantuano_discretization_2007}.

Gao, Brodski, and Mukherjee developed a formulation in which the graph
connection Laplacian is explicitly associated to a flat vector bundle on the
graph and arises from a twisted coboundary operator~\cite{gao_geometry_2016}.
This coboundary operator is not a sheaf coboundary map and has some
difficulties in its definition. These arise from a lack of freedom to choose
the basis for the space of sections over an edge of the graph. Further work by
Gao uses a sheaf Laplacian-like construction to study noninvertible
correspondences between probability distributions on
surfaces~\cite{gao_diffusion_2016}.

Wu et al.~\cite{wu_weighted_2018} have recently proposed a construction they
call a weighted simplicial complex and studied its associated Laplacians. These are
cellular cosheaves where all stalks are equal to a given vector space and
restriction maps are scalar multiples of the identity. Their work discusses the
cohomology and Hodge theory of weighted simplicial complexes, but did not touch
on issues related to the Laplacian spectrum.

% ==========================================================================================
\section{Harmonicity}
\label{sec:harmonicity}
% ==========================================================================================
As a prelude to results about the spectra of sheaf Laplacians, we will discuss
issues related to harmonic cochains on sheaves. While these do not immediately
touch on the spectral properties of the Laplacian, they are closely bound with
its algebraic properties.

% -------------------------------------------------------------------------------------
\subsection{Harmonic Extension}
\label{sec:harmext}
% -------------------------------------------------------------------------------------
%
\begin{proposition}
\label{harmonic_extension}
Let $X$ be a regular cell complex with a weighted cellular sheaf $\Fc$. Let $B
\subseteq X$ be a subcomplex and let $x|_B \in C^k(B;\Fc)$ be an $\Fc$-valued
$k$-cochain specified on $B$. If $H^k(X,B;\Fc) = 0$, then there exists a unique
cochain $x \in C^k(X;\Fc)$ which restricts to $x|_B$ on $B$ and is harmonic on
$S = X \setminus B$.
\end{proposition}
\begin{proof}
A matrix algebraic formulation suffices. Representing $\Delta^k_\Fc$ in block
form as partitioned by $B$ and $S$, the relevant equation is
\[
\begin{bmatrix}
  \Delta^k_\Fc(S,S) & \Delta^k_\Fc(S,B) \\
  \Delta^k_\Fc(B,S) & \Delta^k_\Fc(B,B)
\end{bmatrix}
\begin{bmatrix}
    x|_S \\ x|_B
\end{bmatrix}
=
\begin{bmatrix}
  0 \\ y
\end{bmatrix}.
\]
Since $y$ is indeterminate, we can ignore the second row of the matrix, giving
the equation $\Delta^k_\Fc(S,S)x|_S + \Delta^k_\Fc(S,B)x|_B = 0$. We can write
$\Delta_\Fc^k(S,S) = (\delta^k|_S)^*\delta^k|_S +
((\delta^{k-1})^*|_S)^*(\delta^{k-1})^*|_S$, which is very close to the $k$-th
Hodge Laplacian of the relative cochain complex
\[
  \cdots \to C^{k-1}(X,B;\Fc) \to C^k(X,B;\Fc) \to C^{k+1}(X,B;\Fc) \to \cdots.
\]
Indeed, we can exploit the fact that this is a subcomplex of $C^\bullet(X;\Fc)$
to compute its Hodge Laplacian in terms of the coboundary maps of
$C^\bullet(X;\Fc)$. The coboundary map $\delta_S$ of $C^\bullet(X,B;\Fc)$ is
simply the restriction of the coboundary map $\delta$ of $C^\bullet(X;\Fc)$ to
the subcomplex: $\delta_S^k = \pi_S^{k+1} \delta^k i_S^k$, where $\pi_S^k$ is
the orthogonal projection $C^k(X;\Fc) \to C^k(X,B;\Fc)$ and $i_S^k$ the
inclusion $C^k(X,B;\Fc) \to C^k(X;\Fc)$.
Note that $\pi_S^{k}$ and $i_S^k$ are adjoints, and that $i_S^k \pi_S^k$ is the identity on $\im \delta^{k-1}_S$. 
We may therefore write the Hodge Laplacian of the relative complex as
\begin{multline*}
    \Delta^k(X,B;\Fc) = (\delta^k_S)^*\delta^k_S + \delta^{k-1}_S(\delta^{k-1}_S)^* \\
= \pi_S^k(\delta^k)^*i_S^{k+1} \pi_S^{k+1}\delta^ki_S^k + \pi_S^k \delta^{k-1}i_S^{k-1}\pi_S^{k-1}(\delta^{k-1})^*i_S^k \\
=  \pi_S^k(\delta^k)^*\delta^ki_S^k + \pi_S^k \delta^{k-1}i_S^{k-1}\pi_S^{k-1}(\delta^{k-1})^*i_S^k
.
\end{multline*}
Meanwhile, we can write the submatrix
\[
    \Delta_\Fc^k(S,S) = \pi_S^k(\delta^k)^*\delta^ki_S^k + \pi_S^k \delta^{k-1}(\delta^{k-1})^*i_S^k.
\]
It is then immediate that $\ker(\Delta_\Fc^k(S,S)) \subseteq \ker
\Delta^k(X,B;\Fc)$, so that $\Delta_\Fc^k(S,S)$ is invertible if $H^k(X,B;\Fc)
= 0$.
\end{proof}

If we restrict to up- or down-Laplacians, a harmonic extension always exists,
even if it is not unique. This is because, for instance, $\img
(\delta^k|_S)^*\delta^k|_B \subseteq \img (\delta^k|_S)^*\delta^k|_S$. In
particular, this implies that harmonic extension is always possible for
$0$-cochains, with uniqueness if and only if $H^0(X,B;\Fc) = 0$.

% -------------------------------------------------------------------------------------
\subsection{Kron Reduction}
\label{sec:kron!}
% -------------------------------------------------------------------------------------
\style{Kron reduction} is one of many names given to a process of simplifying
graphs with respect to the properties of their Laplacian on a boundary. If $G$
is a connected graph with a distinguished set of vertices $B$, which we
consider as a sort of boundary of $G$, Proposition~\ref{harmonic_extension}
shows that there is a harmonic extension map $E: \R^B \to \R^{V(G)}$. It is
then possible to construct a graph $G'$ on $B$ such that for every function $x$
on the vertices of $G'$, we have $L_{G'} x = \pi_B L_G E(x)$, where $\pi_B$ is
the orthogonal projection map $\R^{V(G)} \to \R^B$. Indeed, letting $S = V(G)
\setminus B$, we have $E(x)|_S = -L_G(S,S)^{-1} L_G(S,B)x$, so
\[
  L_{G'}x = \pi_B L_G E(x) = L_G(B,B)x -L_G(B,S)L_G(S,S)^{-1} L_G(S,B)x.
\]
Therefore,
\[
  L_{G'} = L_G(B,B)-L_G(B,S)L_G(S,S)^{-1} L_G(S,B),
\]
that is, $L_{G'}$ is the Schur complement of the $(B,B)$ block of $L_G$. It is also the Laplacian of a graph on $B$:

\begin{theorem}[see~\cite{dorfler_kron_2013}]
If $L_G$ is the Laplacian of a connected graph $G$, and $B$ a subset of
vertices of $G$, then $L_{G'} = L_G(B,B)-L_G(B,S)L_G(S,S)^{-1} L_G(S,B)$ is the
Laplacian of a graph with vertex set $B$.
\end{theorem}

A physically-inspired way to understand this result (and a major use of Kron
reduction in practice) is to view it as reducing a network of resistors given
by $G$ to a smaller network with node set $B$ that has the same electrical
behavior on $B$ as the original network. In this guise, Kron reduction is a
high-powered version of the $Y$-$\Delta$ and star-mesh transforms familiar from
circuit analysis. Further discussion of Kron reduction and its various
implications and applications may be found in~\cite{dorfler_kron_2013}.

Can we perform Kron reduction on sheaves? That is, given a sheaf $\Fc$ on a
graph $G$ with a prescribed boundary $B$, can we find a sheaf $\Fc_B$ on a
graph with vertex set $B$ only such that for every $x \in C^0(B;\Fc_B)$ we have
$L_{\Fc_B} x = \pi_{C^0(B;\Fc_B)} L_\Fc E(x)$, where $E(x)$ is the harmonic
extension of $x$ to $G$?

The answer is, in general, no. Suppose we want to remove the vertex $v$ from
our graph, \ie, $B = G \setminus \{v\}$. Let $D_v = \sum_{v \face e}
\Fc_{v\face e}^*\Fc_{v\face e} = L_{v,v}$. To eliminate the vertex $v$ we apply
the condition $(L_{\Fc} (x,E(x)))(v) = 0$, and take a Schur complement,
replacing $L(B,B)$ with $L(B,B) - L(B,v)D_v^{-1}L(v,B)$. This means that we add
to the entry $L(w,w')$ the map $\Fc_{w\face e}^*\Fc_{v\face
e}D_v^{-1}\Fc_{v\face e'}^*\Fc_{w'\face e'}$, where $e$ is the edge between $v$
and $w$, and $e'$ the edge between $v$ and $w'$. This does not in general
translate to a change in the restriction maps for the edge between $w$ and
$w'$. In general, Kron reduction is not possible for sheaves.

\begin{figure}
    {
\centering
\includegraphics[width=3in]{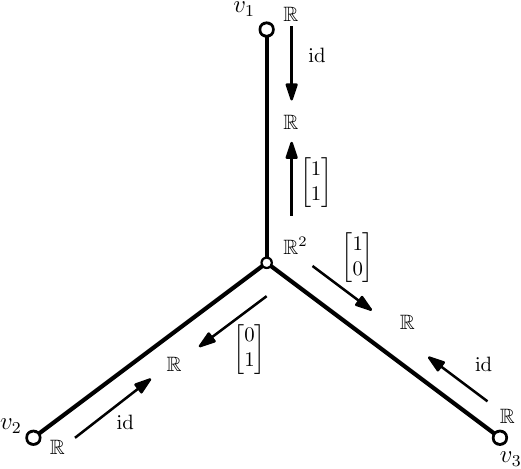}
\caption{A sheaf illustrating the general impossibility of Kron reduction.}\label{fig:nokron}
}
\end{figure}

In particular, if $x\in C^0(G;\Fc)$ is a section of $\Fc$, its restriction to
$B$ must be a section of $\Fc_B$. Conversely, if $x$ is not a section, its
restriction to $B$ cannot be a section of $\Fc_B$. But we can construct sheaves
with a space of sections on the boundary that cannot be replicated with a sheaf
on the boundary vertices only. For instance, take the star graph with three
boundary vertices, with stalks $\R$ over boundary vertices and edges, and
$\R^2$ over the internal vertex. Take as the restriction maps from the central
vertex restriction onto the first and second components, and addition of the
two components. See Figure~\ref{fig:nokron} for an illustration.

Note that a global section of this sheaf is determined by its value on the
central vertex. If we label the boundary vertices counterclockwise starting at
the top, the space of global sections for $\Fc_B$ must have as a basis
\[x_1 = \begin{bmatrix} 1 \\ 1 \\ 0 \end{bmatrix}, \quad x_2 = \begin{bmatrix} 1 \\ 0 \\1 \end{bmatrix}.\]
But there is no sheaf on a graph with vertex set $B$ which has this space of
global sections. To see this, note that if $x_1$ is a section, the map from
$\Fc(v_1)$ to $\Fc(v_3)$ must be the zero map, and similarly for the map from
$\Fc(v_2)$ to $\Fc(v_3)$. Similarly, if $x_2$ is a section, the maps $\Fc(v_1)
\to \Fc(v_2)$ and $\Fc(v_3) \to \Fc(v_2)$ must be zero. But this already shows
that the vector $\begin{bmatrix} 1 & 0 & 0 \end{bmatrix}^T$ must be a section,
giving $\Fc_B$ a three-dimensional space of sections. The problem is that the
internal node allows for constraints between boundary nodes that cannot be
expressed by purely pairwise interactions. This fact is a fundamental
obstruction to Kron reduction for general sheaves.

However, there is a sheaf Kron reduction for sheaves with vertex stalks of
dimension at most 1. This follows from the identification of the Laplacians of
such sheaves as the matrices of factor width two in \S\ref{sec:ident}.

\begin{theorem}
The class of matrices of factor width at most two is closed under taking Schur complements.
\end{theorem}
\begin{proof}
By Theorems 8 and 9 of~\cite{boman_factor_2005}, a matrix $L$ has factor width
at most two if and only if it is symmetric and generalized weakly diagonally
dominant with nonnegative diagonal, that is, there exists a positive diagonal
matrix $D$ such that $DLD$ is weakly diagonally dominant. Equivalently, these
are the symmetric positive semidefinite generalized weakly diagonally dominant
matrices. The class of generalized weakly diagonally dominant matrices
coincides with the class of $H$-matrices, which are shown to be closed under
Schur complements in~\cite{johnson_closure_2005}. Similarly, the class of
symmetric positive definite matrices is closed under Schur complements, so the
intersection of the two classes is also closed.
\end{proof}

% -------------------------------------------------------------------------------------
\subsection{Maximum Modulus Theorem}
\label{sec:mmt}
% -------------------------------------------------------------------------------------
Harmonic $0$-cochains of an $O(n)$-bundle satisfy a local averaging property
which leads directly to a maximum modulus principle.

\begin{lemma}
    Let $G$ be a graph with an $O(n)$-bundle $\Fc$, with constant vertex
    weights $\alpha_v = 1$ and arbitrary edge weights $\alpha_e$ (as defined in
    \S\ref{sec:discvect}). If $x \in C^0(G;\Fc)$ is harmonic at a vertex $v$,
    then
\[x_v = \frac{1}{d_v} \sum_{\substack{v,w \face e\\ v \neq w}} \Fc_{v\face e}^*\Fc_{w\face e}x_w,\] where $d_v = \sum_{v\face e} \norm{\Fc_{v\face e}}^2 = \sum_{v \face e} \alpha_e^2$.
\end{lemma}
\begin{proof}
The block row of $L_\Fc$ corresponding to $v$ has entries $-\Fc_{v\face
e}^*\Fc_{w\face e}$ off the diagonal and $\sum_{v\face e} \Fc_{v\face
e}^*\Fc_{v\face e} = \sum_{v \face e} \norm{\Fc_{v\face e}}^2 \id_{\Fc(v)}$ on
the diagonal. The harmonicity condition is then
\[d_v x_v - \sum_{\substack{v,w \face e\\ v \neq w}} \Fc_{v\face e}^*\Fc_{w\face e} x_w = 0.\]
\end{proof}

\begin{theorem}[Maximum modulus principle]
\label{thm:mmt}
Let $G$ be a graph, and $B$ be a \style{thin} subset of vertices of $G$; that
is, $G \setminus B$ is connected, and every vertex in $B$ is connected to a
vertex not in $B$. Let $\Fc$ be an $O(n)$-bundle on $G$ with $a_v = 1$ for all
$v\in G$, and suppose $x \in C^0(G;\Fc)$ is harmonic on $G\setminus B$. Then if
$x$ attains its maximum stalkwise norm at a vertex in $G \setminus B$, it has
constant stalkwise norm.
\end{theorem}
\begin{proof}
    Note that for a given edge $e = v \sim w$, $\Fc_{v\face e}$ and $\Fc_{u
    \face e}$ are both $\alpha_e$ times an orthogonal map, so $\Fc_{v\face
e}^*\Fc_{w\face e}$ is $\alpha_e^2$ times an orthogonal map. 
Let $v \in G \setminus B$ and suppose $\norm{x_v} \geq \norm{x_w}$ for all $w
\in G$. Then this holds in particular for neighbors of $v$, so that we have
\begin{multline*}
\norm{x_v} = \frac{1}{d_v} \norm{\sum_{\substack{v,w \face e\\ v \neq w}} \Fc_{v\face e}^*\Fc_{w\face e}x_w} \leq \frac{1}{d_v} \sum_{\substack{v,w \face e\\ v \neq w}} \norm{\Fc_{v\face e}^*\Fc_{w\face e}x_w} \\
= \frac{1}{d_v} \sum_{\substack{v,w \face e\\ v \neq w}} \alpha_e^2 \norm{x_w} \leq \frac{1}{d_v} \sum_{v\face e} \alpha_e^2 \norm{x_v} = \norm{x_v},
\end{multline*}
Equality holds throughout, which, combined with the assumption that $\norm{x_v}
\geq \norm{x_w}$ for all $w$, forces $\norm{x_v} = \norm{x_w}$ for $w \sim v$.
We then apply the same argument to every vertex in $G\setminus B$ adjacent to
$v$, and, after iterating, the region of constant stalkwise norm extends to all
of $G \setminus B$ because this subgraph is connected. But since every vertex
$b \in B$ is adjacent to some vertex $w \in G \setminus B$, the same argument
applied to the neighborhood of $w$ forces $\norm{x_b} = \norm{x_w}$. So any
harmonic function that attains its maximum modulus on $G \setminus B$ has
constant modulus.
\end{proof}

\begin{corollary}
    Let $B$ be a thin subset of vertices of $G$, and $\Fc$ an $O(n)$-bundle on
    $G$ as before. If $x \in C^0(G;\Fc)$ is harmonic on $G\setminus B$, then it
    attains its maximum modulus on $B$.
\end{corollary}

The constant sheaf on a graph is an $O(n)$-bundle, so this result gives a
maximum modulus principle for harmonic functions on the vertices of a graph. A
slightly stronger result in this vein, involving maxima and minima of $x$, is
discussed in~\cite{sunada_discrete_2008}. The thinness condition for $B$ is not
strictly necessary for the corollary to hold --- there are a number of
potential weakenings of the condition. For instance, we might simply require
that there exists some $w \in B$ such that for every vertex $v \in G \setminus
B$ there exists a path from $v$ to $w$ not passing through $B$.

% ==========================================================================================
\section{Spectra of Sheaf Laplacians}
\label{sec:spectralbasics}
% ==========================================================================================
The results in this section are straightforward generalizations and extensions
of familiar results from spectral graph theory. Most are not particularly
difficult, but they illustrate the potential for lifting nontrivial notions
from graphs and complexes to sheaves.

It is useful to note a few basic facts about the spectra of Laplacians arising from Hodge theory.

\begin{proposition}
The nonzero spectrum of $\Delta^k$ is the disjoint union of the nonzero spectra of $\Delta^k_+$ and $\Delta^k_-$.
\end{proposition}
\begin{proof}
We take advantage of the Hodge decomposition, noting that $C^k(X;\Fc) = \ker
\Delta^k \oplus \im \Delta^k_- \oplus \im \Delta^k_+$. This is an orthogonal
decomposition, and $ = 0$ as well as $\Delta^k_-|_{(\im \Delta^k_+)} = 0$.
Further, since $\ker \Delta^k = \ker \Delta_+ \cap \ker \Delta_-$, both
restrict to zero on the kernel of $\Delta^k$. We therefore see that $\Delta^k$
is the orthogonal direct sum $0|_{\ker \Delta^k}\oplus \Delta^k_+|_{(\im
\Delta^k_+)} \oplus \Delta^k_-|_{(\im \Delta^k_-)}$, and hence the spectrum of
$\Delta^k$ is the union of the spectra of these three operators.
\end{proof}

\begin{proposition}\label{prop:updownspectraequal}
The nonzero eigenvalues of $\Delta^k_+$ and $\Delta^{k+1}_-$ are the same.
\end{proposition}
\begin{proof}
We have $\Delta^k_+ = (\delta^k)^*\delta^k$ and $\Delta^{k+1}_- =
\delta^k(\delta^k)^*$. The eigendecompositions of these matrices are determined
by the singular value decomposition of $\delta^k$, and the nonzero eigenvalues
are precisely the squares of the nonzero singular values of $\delta^k$.

\end{proof}

One reason for the study of the normalized graph Laplacian is that its spectrum
is bounded above by 2~\cite{chung_spectral_1992}, and hence normalized
Laplacian spectra of different graphs can be easily compared. A similar result holds for
up-Laplacians of normalized simplicial complexes~\cite{horak_spectra_2013}: the
eigenvalues of the degree-$k$ up-Laplacian of a normalized simplicial complex
are bounded above by $k+2$. This fact extends to normalized sheaves on
simplicial complexes. This result and others in this paper will rely on the
Courant-Fischer theorem, which we state here for reference.
\begin{definition}
    Let $A$ be a self-adjoint operator on a Hilbert space $\Vspace$. If $x \in
    \Vspace$, the \style{Rayleigh quotient} corresponding to $x$ and $A$ is
    \[
        R_A(x) = \frac{\ip{x,Ax}}{\ip{x,x}}.
    \]
\end{definition}

\begin{theorem}[Courant-Fischer]
    Let $A$ be an $n \times n$ Hermitian matrix with eigenvalues $\lambda_1 \leq
    \lambda_2 \leq \cdots \leq \lambda_n$. Then
    \[\lambda_k = \min_{\dim V = k} \max_{x \in
        V} R_A(x) = \max_{\dim V = n-k+1}
    \min_{x \in V} R_A(x).\]
\end{theorem}
The proof is immediate once one uses the fact that $A$ is unitarily equivalent
to a diagonal matrix.

\begin{proposition}
Suppose $\Fc$ is a normalized sheaf on a simplicial complex $X$. The
eigenvalues of the degree $k$ up-Laplacian $L^k_\Fc$ are bounded above by
$k+2$.
\end{proposition}
\begin{proof}
By the Courant-Fischer theorem, the largest eigenvalue of $L^k_\Fc$ is equal to
\begin{multline*}
\max_{x \in C^k(X;\Fc)} \frac{\ip{x,L^k_\Fc x}}{\ip{x,x}}
=
\max_{x \perp \ker \delta^k} \frac{\ip{\delta^kx,\delta^k x}}{\sum\limits_{\dim \sigma = k} \ip{\delta^k x_\sigma,\delta^k x_\sigma}}
\\
=
\max_{x \perp \ker \delta^k} \frac{\sum\limits_{\dim \tau = k+1} \sum\limits_{\sigma,\sigma'\face\tau} [\sigma:\tau][\sigma':\tau]\ip{\Fc_{\sigma\face\tau} x_\sigma , \Fc_{\sigma'\face\tau} x_{\sigma'}}}{\sum\limits_{\dim \sigma = k} \sum\limits_{\sigma\face\tau} \ip{\Fc_{\sigma\face\tau}x_\sigma,\Fc_{\sigma\face\tau}x_\sigma}}
.
\end{multline*}
Note that for $\sigma \neq \sigma'$,
\begin{multline*}
  [\sigma:\tau][\sigma':\tau]\ip{\Fc_{\sigma\face\tau} x_\sigma , \Fc_{\sigma'\face\tau} x_{\sigma'}}
  \leq
  \norm{\Fc_{\sigma\face\tau}x_\sigma}\norm{\Fc_{\sigma'\face\tau}x_{\sigma'}}
  \\
  \leq
  \frac{1}{2}\left(\norm{\Fc_{\sigma\face\tau}x_\sigma}^2
  +
  \norm{\Fc_{\sigma'\face\tau}x_{\sigma'}}^2\right)
\end{multline*}
by the Cauchy-Schwarz inequality. In particular, then, the term of the
numerator corresponding to each $\tau$ of dimension $k+1$ is bounded above by
\[
  \sum_{\sigma\face\tau} \norm{\Fc_{\sigma\face\tau}x_\sigma}^2 +  \frac{1}{2}
  \sum_{\sigma\neq \sigma'\face\tau} \left(
  \norm{\Fc_{\sigma\face\tau}x_\sigma}^2 +
  \norm{\Fc_{\sigma'\face\tau}x_{\sigma'}}^2\right) =  (k+2) \sum_{\sigma\face\tau} \norm{\Fc_{\sigma\face\tau} x_\sigma}^2
  ,
\]
by counting the number of times each term $\norm{\Fc_{\sigma \face
\tau}x_\sigma}^2$ appears in the sum.
Meanwhile, the denominator is equal to $\displaystyle\sum_{\dim \tau = k+1}
\sum_{\sigma\face\tau} \norm{\Fc_{\sigma\face\tau}x_\sigma}^2$, so the Rayleigh
quotient is bounded above by $k+2$.
\end{proof}

% -------------------------------------------------------------------------------------
\subsection{Eigenvalue Interlacing}
\label{sec:interlacing}
% -------------------------------------------------------------------------------------

\begin{definition}
Let $A$, $B$ be $n \times n$ matrices with real spectra. Let $\lambda_1\leq
\lambda_2 \leq \cdots \leq \lambda_n$ be the eigenvalues of $A$ and $\mu_1 \leq
\mu_2 \leq \cdots \leq \mu_n$ be the eigenvalues of $B$. We say the eigenvalues
of $A$ are $\mathbf{(p,q)}$-\style{interlaced} with the eigenvalues of $B$ if
for all $k$, $\lambda_{k-p} \leq \mu_k \leq \lambda_{k+q}$. (We let
$\lambda_{k} = \lambda_1$ for $k < 1$ and $\lambda_k = \lambda_n$ for $k > n$.)
\end{definition}

The eigenvalues of low-rank perturbations of symmetric positive semidefinite
matrices are related by interlacing. The following is a standard result:

\begin{theorem}
Let $A$ and $B$ be positive semidefinite matrices, with $\rank B = t$. Then the
eigenvalues of $A$ are $(t,0)$-interlaced with the eigenvalues of $A - B$.
\end{theorem}
\begin{proof}
Let $\mu_k$ be the $k$-th largest eigenvalue of $A-B$ and $\lambda_k$ the $k$-th largest eigenvalue of $A$.
By the Courant-Fischer theorem, we have
\begin{align*}
\mu_k
&= \min_{\dim Y = k} \left(\max_{y \in Y,y\neq 0} \frac{\ip{y,A y} - \ip{y, By}}{\ip{y, y}}\right)\\
&\geq \min_{\dim Y = k} \left(\max_{y \in Y \cap \ker B,y\neq 0} \frac{\ip{y,A y}}{\ip{y, y}}\right)\\
&\geq \min_{\dim Y = k - t } \left(\max_{y \in Y,y\neq 0} \frac{\ip{y,A y}}{\ip{y, y}}\right) = \lambda_{k-t}
\end{align*}
and
\begin{align*}
\lambda_k
& = \min_{\dim Y = k} \left(\max_{y \in Y,y\neq 0} \frac{\ip{y,A y}}{\ip{y, y}}\right) \\
& \geq \min_{\dim Y = k} \left(\max_{y \in Y,y\neq 0} \frac{\ip{y,A y} - \ip{y, By}}{\ip{y, y}}\right) = \mu_{k}.
\end{align*}
\end{proof}

This result is immediately applicable to the spectra of sheaf Laplacians under
the deletion of cells from their underlying complexes. The key part is the
interpretation of the difference of the two Laplacians as the Laplacian of a
third sheaf.\footnote{Such subtle moves are part and parcel of a
sheaf-theoretic perspective.} Let $\Fc$ be a sheaf on $X$, and let $C$ be an
upward-closed set of cells of $X$, with $Y = X \setminus C$. The inclusion map
$i: Y \to X$ induces a restriction of $\Fc$ onto $Y$, the pullback sheaf
$i^*\Fc$. Consider the Hodge Laplacians $\Delta^k_\Fc$ and $\Delta^k_{i^*\Fc}$.
If $C$ contains cells of dimension $k$, these matrices are different sizes, but
we can derive a relationship by padding $\Delta^k_{i^*\Fc}$ with zeroes.
Equivalently, this is the degree-$k$ Laplacian of $\Fc$ with the restriction
maps incident to cells in $C$ set to zero.

\begin{proposition}
Let $\Gc$ be the sheaf on $X$ with the same stalks as $\Fc$ but with all
restriction maps between cells not in $C$ set to zero.
The eigenvalues of $\Delta^k_{i^*\Fc}$ are $(t,0)$-interlaced with the
eigenvalues of $\Delta^k_\Fc$, where $t = \codim H^k(X;\Gc) =  \dim C^k(X; \Fc)
- \dim H^k(X;\Gc)$.
\end{proposition}

Similar results can be derived for the up- and down-Laplacians. Specializing to
graphs, interlacing is also possible for the normalized degree 0 sheaf
Laplacian. The Rayleigh quotient for the normalized Laplacian $\mathcal
L_{i^*\Fc}$ is
\[
  \frac{\ip{x, D_{i^*\Fc}^{-1/2}L_{i^*\Fc}D_{i^*\Fc}^{-1/2}x}}{\ip{x, x}}
  =
  \frac{\ip{y,L_{i^*\Fc} y}}{\ip{y, D_{i^*\Fc}y}}
  =
  \frac{\ip{y,L_\Fc y} - \ip{y, L_\Gc y}}{\ip{y, D_\Fc y} - \ip{y, D_\Gc y}},
\]
where we let $y = D_{i^*\Fc}^{-1/2}x$. Then if $\{\lambda_k\}$ are the ordered eigenvalues of $\mathcal L_{\Fc}$ and $\{\mu_k\}$ are the ordered eigenvalues of $\mathcal L_{i^*\Fc}$, we have
\begin{align*}
\mu_k
&= \min_{\dim Y = k} \left(\max_{y \in Y,y\neq 0} \frac{\ip{y,L_\Fc y} - \ip{y, L_\Gc y}}{\ip{y, D_\Fc y} - \ip{y, D_\Gc y}}\right)\\
&\geq \min_{\dim Y = k} \left(\max_{y \in Y \cap H^0(X;\mathcal G),y\neq 0} \frac{\ip{y,L_\Fc y}}{\ip{y, D_\Fc y} - \ip{y, D_\Gc y}}\right)\\
&\geq \min_{\dim Y = k} \left(\max_{y \in Y \cap H^0(X;\mathcal G),y\neq 0} \frac{\ip{y,L_\Fc y}}{\ip{y, D_\Fc y}}\right)\\
&\geq \min_{\dim Y = k - t } \left(\max_{y \in Y,y\neq 0} \frac{\ip{y,L_\Fc y}}{\ip{y, D_\Fc y}}\right) = \lambda_{k-t}
\end{align*}

\begin{align*}
\mu_k
&= \max_{\dim Y = n-k+1} \left(\min_{y \in Y,y\neq 0} \frac{\ip{y,L_\Fc y} -
\ip{y, L_\Gc y}}{\ip{y, D_\Fc y} - \ip{y, D_\Gc y}}\right)\\
&\leq \max_{\dim Y = n-k+1} \left(\min_{y \in Y \cap H^0(X;\mathcal G),y\neq 0} \frac{\ip{y,L_\Fc y}}{\ip{y, D_\Fc y} - \ip{y, D_\Gc y}}\right)\\
&\leq \max_{\dim Y = n-k+1} \left(\min_{y \in Y \cap H^0(X;\mathcal G),y\neq 0} \frac{\ip{y,L_\Fc y}}{\ip{y, D_\Fc y}}\right)\\
&\leq \max_{\dim Y = n-k-t+1} \left(\min_{y \in Y,y\neq 0} \frac{\ip{y,L_\Fc y}}{\ip{y, D_\Fc y}}\right) = \lambda_{k+t}
\end{align*}
Therefore, the eigenvalues of the normalized Laplacians are $(t,t)$-interlaced.
This generalizes interlacing results for normalized graph Laplacians.

% -------------------------------------------------------------------------------------
\subsection{Sheaf Morphisms}
\label{sec:shmorph}
% -------------------------------------------------------------------------------------

\begin{proposition}
Suppose $\varphi: \Fc \to \Gc$ is a morphism of weighted sheaves on a regular
cell complex $X$. If $\varphi^{k+1}$ is a unitary map, then $L^k_\Fc =
(\varphi^k)^*L^k_\Gc \varphi^k$.
\end{proposition}
\begin{proof}
    The commutativity condition $\varphi^{k+1} \delta^\Fc = \delta^\Gc
    \varphi^k$ implies that $(\delta^\Fc)^* (\varphi^{k+1})^* \varphi^{k+1}
    \delta^\Fc = (\varphi^k)^*(\delta^\Gc)^*\delta^\Gc \varphi^k =
    (\varphi^k)^*L^k_\Gc \varphi^k$. Thus if $(\varphi^{k+1})^* \varphi^{k+1}=
    \id_{C^{k+1}(X;\Fc)}$, we have $L^k_\Fc = (\varphi^k)^*L^k_\Gc \varphi^k$.
    This condition holds if $\varphi^{k+1}$ is unitary.
\end{proof}
An analogous result holds for the down-Laplacians of $\Fc$, and these combine
to a result for the full Hodge Laplacians.

% -------------------------------------------------------------------------------------
\subsection{Cell Complex Morphisms}
\label{sec:celmorph}
% -------------------------------------------------------------------------------------
The following constructions are restricted to \style{locally injective}
cellular morphisms, as discussed in \S\ref{sec:cellcomp}. Recall that
under these morphisms, cells map to cells of the same dimension and the
preimage of the star of a cell is a disjoint union of subcomplexes, on each of
which the map acts injectively. The sheaf Laplacian is invariant with
respect to pushforwards over such maps:

\begin{proposition}
\label{prop:pushforward}
Let $X$ and $Y$ be cell complexes, and let $f: X \to Y$ be a locally injective
cellular morphism. If $\Fc$ is a sheaf on $X$, the $k$th coboundary Laplacian
corresponding to $f_* \Fc$ on $Y$ is the same (up to a unitary change of basis)
as the $k$th coboundary Laplacian of $\Fc$ on $X$.
\end{proposition}

\begin{corollary}
The sheaves $\Fc$ and $f_*\Fc$ are isospectral for the coboundary Laplacian.
\end{corollary}
\begin{proof}
There is a canonical isometry $f_k: C^k(X,\Fc)\to C^k(Y,f_*\Fc)$, which is
given on stalks by the obvious inclusion $f_\sigma: \Fc(\sigma) \to
f_*\Fc(f(\sigma)) = \bigoplus_{f(\tau) = f(\sigma)} \Fc(\tau)$. For
$\sigma\face\sigma'$, $f_\sigma$ commutes with the restriction map
$\Fc_{\sigma\face\sigma'}$ and hence $f_k$ commutes with the coboundary map.
But this implies that:
\[
  L_{f_*\Fc}^k
  = (\delta_{f_*\Fc}^k)^* \delta_{f_*\Fc}^k
  = (\delta_{f_*\Fc}^k)^* f_{(k+1)}^* f_{(k+1)} \delta_{f_*\Fc}^k
  = f_k^* (\delta_\Fc^k)^* \delta_\Fc^k f_k
  = f_k^* L_\Fc^k f_k
  .
\]
\end{proof}
General locally injective maps behave nicely with sheaf pushforwards, and
covering maps behave well with sheaf pullbacks. Recall that
a covering map of cell complexes is a locally injective map $f:C \to X$ such that for every cell $\sigma \in X$, $f$ is an isomorphism on the disjoint components of $f^{-1}(\st(\sigma))$.

\begin{proposition}
Let $f: C \to X$ be a covering map of cell complexes, with $\Fc$ a sheaf on
$X$. Then for any $k$, the spectrum of $L_\Fc^k$ is contained in the spectrum
of $L_{f^*\Fc}^k$.
\end{proposition}
\begin{proof}
Consider the lifting map $\Lift: C^k(X;\Fc) \to C^k(C;f^*\Fc)$ given by $x
\mapsto x \circ f_k$. This map commutes with $\delta$ and $\delta^*$. The
commutativity with $\delta$ follows immediately from the proof of the
contravariant functoriality of cochains. The commutativity with $\delta^*$ is
more subtle, and relies on the fact that $f$ is a covering map.

For $y \in C^k(C;f^*\Fc)$ and $x \in C^{k+1}(X;\Fc)$, we have
\begin{align*}
    \ip{y,\delta^*\Lift x} = \ip{\delta y, \Lift x} &= \!\!\sum_{\substack{\sigma',\tau' \in P_C \\ \sigma'\face\tau'}} \![\sigma':\tau'] \,\ip{f^*\Fc_{\sigma'\face\tau'}(y_{\sigma'}), (\Lift x)_{\tau'}} \\
                                                    &=\! \sum_{\substack{\sigma,\tau \in P_X \\ \sigma\face\tau}} \![\sigma: \tau] \!\!\sum_{\sigma' \in f^{-1}(\sigma)} \!\ip{\Fc_{\sigma\face\tau}(y_{\sigma'}),x_\tau} \\
                                                    &=
                                                    \sum_{\substack{\sigma,\tau
                                                    \in P_X \\
                                            \sigma\face\tau}} [\sigma :\tau]
                                            \ip{\Fc_{\sigma\face\tau} (\Lift^*
                                            y)_\sigma,x_\tau} \\
                                                    &= \ip{\delta
                                        \Lift^* y, x} = \ip{y,\Lift \delta^* x}.
\end{align*}

Now, if $L^k_\Fc x = \lambda x$, we have $L^k_{f^*\Fc}\Lift x =
(\delta_{f^*\Fc}^k)^*\delta_{f^*\Fc}^k \Lift x =
\Lift(\delta_\Fc^k)^*\delta_\Fc^k x = \Lift L^k_\Fc x = \lambda \Lift x$, so
$\lambda$ is an eigenvalue of $L^k_{f^*\Fc}$.
\end{proof}

Even if $f: Y \to X$ is not quite a covering map, it is still possible to get
some information about the spectrum of $f^*\Fc$. For instance, for
dimension-preserving cell maps with uniform fiber size we have a bound on the
smallest nontrivial eigenvalue of the pullback:

\begin{proposition}
Suppose $f: Y \to X$ is a \style{dimension-preserving} map of regular cell
complexes such that for $\dim(\sigma) = d$, $\abs{f^{-1}(\sigma)} = \ell_d$ is
constant, and let $\Fc$ be a sheaf on $X$. If $\lambda_k(\Fc)$ is the smallest
nontrivial eigenvalue of $L_{\Fc}^d$, then $\lambda_k({\Fc}) \geq
\frac{\ell_d}{\ell_{d+1}}\lambda_k({f^* \Fc})$.
\end{proposition}
\begin{proof}
Let $x$ be an eigenvector corresponding to $\lambda_k({\Fc})$. Note that since
every fiber is the same size, the lift $\Lift$ preserves the inner product up
to a scaling. That is, if $y$ and $z$ are  $d$-cochains, $\ip{\Lift y, \Lift z}
= \ell_d \ip{y,z}$. This means that the pullback of $x$ is orthogonal to the
pullback of any cochain in the kernel of $L_\Fc$. Therefore, we have
\[
\lambda_k(\Fc) 
= \frac{\ip{\delta x,\delta x}}{\ip{x,x}} 
= \frac{\ell_{d}\ip{\Lift \delta x, \Lift \delta x}}{\ell_{d+1} \ip{\Lift x, \Lift x}} 
= \frac{\ell_{d}\ip{\delta \Lift x, \delta \Lift x}}{\ell_{d+1}\ip{\Lift x, \Lift x}} 
\geq \frac{\ell_{d}}{\ell_{d+1}}\lambda_k({f^*\Fc}).
\]
\end{proof}

% -------------------------------------------------------------------------------------
\subsection{Product Complexes}
\label{sec:procom}
% -------------------------------------------------------------------------------------
If $X$ and $Y$ are cell complexes, their product $X \times Y$ is a cell complex
with cells $\sigma \times \tau$ for $\sigma \in X$, $\tau \in Y$, and incidence
relations $(\sigma \times \tau)\face (\sigma' \times \tau')$ whenever
$\sigma\face\sigma'$ and $\tau\face\tau'$. The dimension of $\sigma \times
\tau$ is $\dim(\sigma) + \dim(\tau)$. The complex $X \times Y$ possesses
projection maps $\pi_X$ and $\pi_Y$ onto $X$ and $Y$.

\begin{definition}
If $\Fc$ and $\Gc$ are sheaves on $X$ and $Y$, respectively, their
\style{product} is the sheaf $\Fc \boxtimes \Gc = \pi_X^* \Fc \otimes \pi_Y^*
\Gc$. Equivalently, we have $(\Fc \boxtimes \Gc)(\sigma \times \tau) =
\Fc(\sigma) \otimes \Fc(\tau)$ and $(\Fc \boxtimes \Gc)_{\sigma \times
\tau\face\sigma' \times \tau'} = \Fc_{\sigma\face\sigma'} \otimes
\Gc_{\tau\face\tau'}$.
\end{definition}

\begin{proposition}
If $L_\Fc$ and $L_\Gc$ are the degree-0 Laplacians of $\Fc$ and $\Gc$, the
degree-0 Laplacian of $\Fc \boxtimes \Gc$ is $L_{\Fc \boxtimes \Gc} =
\id_{C^0(X;\Fc)} \otimes L_\Gc + L_\Fc \otimes \id_{C^0(Y;\Gc)}$.
\end{proposition}
\begin{proof}
The vector space $C^1(X\times Y;\Fc\boxtimes \Gc)$ has a natural decomposition
into two subspaces: one generated by stalks of the form $\Fc(v) \otimes \Gc(e)$
for $v$ a vertex of $X$ and $e$ an edge of $Y$, and another generated by stalks
of the opposite form $\Fc(e) \otimes \Gc(v)$. This induces an isomorphism
\[
  C^1(X\times Y;\Fc\boxtimes \Gc) \cong (C^0(X;\Fc) \otimes C^1(Y;\Gc)) \oplus (C^1(X;\Fc) \otimes C^0(Y;\Gc)).
\]
Then the coboundary map of $\Fc \boxtimes \Gc$ can be written as the block matrix
\[
\delta_{\Fc \boxtimes \Gc}
=
\begin{bmatrix}
  \id_{C^0(X;\Fc)} \otimes \delta_\Gc \\
  \delta_\Fc \otimes \id_{C^0(Y;\Gc)}
\end{bmatrix}.
\]
A quick computation then gives $L_{\Fc \boxtimes \Gc} = \delta_{\Fc \boxtimes
\Gc}^*\delta_{\Fc \boxtimes \Gc} = \id_{C^0(X;\Fc)} \otimes
\delta_\Gc^*\delta_\Gc + \delta_\Fc^*\delta_\Fc \otimes \id_{C^0(Y;\Gc)} =
\id_{C^0(X;\Fc)} \otimes L_\Gc + L_\Fc \otimes \id_{C^0(Y;\Gc)}$.
\end{proof}

\begin{corollary}
If the spectrum of $L_\Fc$ is $\{\mu_i\}_i$ and the spectrum of $L_\Gc$ is
$\{\lambda_j\}_j$, then the spectrum of $L_{\Fc \boxtimes \Gc}$ is $\{\mu_i +
\lambda_j\}_{ij}$.
\end{corollary}

For higher degree Laplacians, this relationship becomes more complicated. For
instance, the degree-1 up-Laplacian is computed as follows:

\[
\delta_{\Fc \boxtimes \Gc}^1
=
\begin{bmatrix}
  \id_{C^0(X;\Fc)} \otimes \delta^1_\Gc & 0\\
  \delta^0_\Fc \otimes \id_{C^1(Y;\Gc)} & \id_{C^1(X;\Fc)}\otimes \delta^0_\Gc \\
  0 & \delta^1_\Fc\otimes \id_{C^0(Y;\Gc)}
\end{bmatrix}.
\]
\[
L^1_{\Fc \boxtimes \Gc}
=
\begin{bmatrix}
  \id_{C^0(X;\Fc)} \otimes L_\Gc^1 + L_\Fc^0\otimes \id_{C^1(Y;\Gc)} & (\delta^0_\Fc)^*\otimes \delta^0_\Gc \\
  \delta^0_\Fc\otimes(\delta^0_\Gc)^* & \id_{C^1(X;\Fc)}\otimes L^0_\Gc + L^1_\Fc \otimes \id_{C^0(Y;\Gc)}
\end{bmatrix}
.
\]
Because $L^1_{\Fc \boxtimes \Gc}$ is given by a block matrix in terms of
various Laplacians and coboundary maps of $\Fc$ and $\Gc$, computing its
spectrum is more involved. When $X$ and $Y$ are graphs, this simplifies
significantly and it is possible to compute the spectrum in terms of the
spectra of $\Fc$ and $\Gc$.

\begin{proposition}
Suppose $X$ and $Y$ are graphs, with $\Fc$ and $\Gc$ sheaves on $X$ and $Y$. If
$v_\Fc$ is an eigenvector of $L_\Fc^0$ with eigenvalue $\lambda$ and $v_\Gc$ an
eigenvector of $L_\Gc^0$ with eigenvalue $\mu$, then the vector $v_{\Fc
\boxtimes \Gc} = \begin{bmatrix} \sqrt{\frac{\lambda}{\mu}} v_\Fc \otimes
\delta_\Gc^0 v_\Gc \\ \sqrt{\frac{\mu}{\lambda}} \delta^0_\Fc v_\Fc \otimes
v_\Gc \end{bmatrix}$ is an eigenvector of $L_{\Fc \boxtimes \Gc}^1$ with
eigenvalue $\lambda + \mu$.
\end{proposition}
\begin{proof}
A computation.
\begin{multline*}
L_{\Fc \boxtimes \Gc}^1 v_{\Fc \boxtimes \Gc}
=
\begin{bmatrix}
  L_\Fc^0\otimes \id_{C^1(Y;\Gc)} & (\delta^0_\Fc)^*\otimes \delta^0_\Gc \\
  \delta^0_\Fc\otimes(\delta^0_\Gc)^* & \id_{C^1(X;\Fc)}\otimes L^0_\Gc
\end{bmatrix}
\begin{bmatrix}
  \sqrt{\frac{\lambda}{\mu}} v_\Fc \otimes \delta_\Gc^0 v_\Gc \\
  \sqrt{\frac{\mu}{\lambda}} \delta^0_\Fc v_\Fc \otimes v_\Gc
\end{bmatrix}
\\
=
\begin{bmatrix}
  \left(\sqrt{\frac{\lambda}{\mu}} +\sqrt{\frac{\mu}{\lambda}}\right)\lambda v_\Fc \otimes \delta^0_\Gc v_\Gc \\
  \left(\sqrt{\frac{\lambda}{\mu}} +\sqrt{\frac{\mu}{\lambda}}\right)\mu \delta^0_\Fc v_\Fc \otimes v_\Gc
\end{bmatrix}
=
(\lambda + \mu)
\begin{bmatrix}
  \sqrt{\frac{\lambda}{\mu}} v_\Fc \otimes \delta_\Gc^0 v_\Gc \\
  \sqrt{\frac{\mu}{\lambda}} \delta^0_\Fc v_\Fc \otimes v_\Gc
\end{bmatrix}.
\end{multline*}
\end{proof}

A simpler way to obtain nearly the same result is to recall from
Proposition~\ref{prop:updownspectraequal} that $(\delta^1_{\Fc \boxtimes
\Gc})^*\delta^1_{\Fc \boxtimes \Gc}$ and $\delta^1_{\Fc \boxtimes \Gc}
(\delta^1_{\Fc \boxtimes \Gc})^*$ have the same spectrum up to the multiplicity
of zero. But when $X$ and $Y$ are graphs,
\[
    \delta^1_{\Fc \boxtimes \Gc} (\delta^1_{\Fc \boxtimes \Gc})^* =
    (\Delta^1_-)_\Fc \otimes \id_{C^1(Y;\Gc)} + \id_{C^1(X;\Fc)} \otimes (\Delta^1_-)_\Gc
,
\]
and since the nonzero eigenvalues of $(\Delta^1_-)_\Fc$ are the same as those of
$L^0_\Fc$, we obtain a correspondence between eigenvalues of $L^1_{\Fc \boxtimes
\Gc}$ and sums of eigenvalues of $L^0_\Fc$ and $L^0_\Gc$.

However, for higher-dimensional complexes and higher-degree Laplacians, there
appears to be no general simple formula giving the spectrum of $\Fc \boxtimes
\Gc$ in terms of the spectra of $\Fc$ and $\Gc$.  Indeed, we suspect
that no such formula can exist, \ie, that the spectrum of $\Fc \boxtimes \Gc$
is not in general determined by the spectra of $\Fc$ and $\Gc$.

% ==========================================================================================
\section{Effective Resistance}
\label{sec:resistance}
% ==========================================================================================

Effective resistance is most naturally defined for weighted cosheaves, but
since every weighted sheaf has a canonical dual weighted cosheaf, the following
definition applies to weighted sheaves as well.

\begin{definition}
    Let $\Fc$ be a cosheaf on a cell complex $X$, and let $a, b \in \ker
    (\partial_k)_\Fc$ be homologous $k$-cycles. The \style{effective
    resistance} $R_{\text{eff}}(a,b)$ is given by the solution to the
    optimization problem
\begin{equation}
\min_{c \in C_{k+1}(X;\Fc)}\, \,\norm{c}^2 \quad \text{s.t. } \,\,\partial_{k+1} c = b-a. \label{eqn:eff_res}
\end{equation}
\end{definition}

\begin{proposition}
Cosheaf effective resistance may be computed by the formula
$R_{\text{eff}}(a,b) = \ip{b-a, (L_\Fc^{k})^\dagger (b-a)}$, where $L_\Fc^k$ is
the up Laplacian $\partial_{k+1}\partial_{k+1}^*$.
\end{proposition}
\begin{proof}
By a standard result about matrix pseudoinverses, $(L_\Fc^{k})^\dagger =
(\partial_{k+1}^\dagger)^*\partial_{k+1}^\dagger$, so $\ip{b-a,
(L_\Fc^{k})^\dagger (b-a)} =
\ip{\partial_{k+1}^\dagger(b-a),\partial_{k+1}^\dagger(b-a)} =
\norm{\partial_{k+1}^\dagger(b-a)}^2$. But if $b-a \in \im \partial_{k+1}$,
then $\partial_{k+1}^\dagger(b-a)$ is the minimizer of the optimization problem~\eqref{eqn:eff_res}.
\end{proof}

The condition that $a$ and $b$ be homologous becomes trivial if $H_k(X;\Fc) =
0$. Note that if $\Fc$ is the constant cosheaf on a graph, a $0$-cycle
supported on a vertex $v$ is homologous to a $0$-cycle supported on a vertex
$w$ if their values are the same. This means that the definition of cosheaf
effective resistance recovers the definition of graph effective resistance.

For any cosheaf on a cell complex, we can get a measure of effective resistance
over a $(k+1)$-cell $\sigma$ by using the boundary map restricted to $\sigma$.
Any choice of $c \in \Fc(\sigma)$ gives an equivalence class of pairs of
homologous $k$-cycles supported inside the boundary of $\sigma$. That is, we
can decompose $\partial c$ into the sum of two $k$-cycles in a number of
equivalent ways. For instance, if $\sigma$ is a 1-simplex with two distinct
incident vertices, there is a natural decomposition of $\partial c$ into a sum
of two $0$-cycles, one supported on each vertex. This gives a quadratic form on
$\Fc(\sigma)$: $R_{\text{eff}}(\sigma)(x) = \ip{\partial x,L_{\Fc}^\dagger
\partial x}$. The choice of decomposition does not affect the quadratic form.
Of course, by the inner product pairing, this quadratic form can be represented as a matrix
$(\partial|_{\Fc(\sigma)})^*L_{\Fc}^\dagger \partial|_{\Fc(\sigma)}$. In
particular, this defines a matrix-valued effective resistance over an edge for a
cosheaf on a graph.

% -------------------------------------------------------------------------------------
\subsection{Sparsification}\label{sec:spars}
% -------------------------------------------------------------------------------------
Graph effective resistance has also attracted attention due to its use in graph
sparsification. The goal when sparsifying a graph $G$ is to find a graph $H$
with fewer edges that closely preserves properties of $G$. One important
property we might wish to preserve is the size of boundaries of sets of
vertices. If $S$ is a set of vertices, we let $\abs{\partial S}$ be the sum of
the weights of edges between $S$ and its complement. We can compute this in
terms of the Laplacian of $G$: $\abs{\partial S} =
\mathbf{1}_S^TL_G\mathbf{1}_S$, where $\mathbf{1}_S$ is the indicator vector on
the set $S$. If $H$ well approximates $G$ in this sense, we
would like to have a relation like $(1-\epsilon)\mathbf{1}_S^TL_G\mathbf{1}_S
\leq \mathbf{1}_S^TL_H\mathbf{1}_S \leq (1+\epsilon)
\mathbf{1}_S^TL_G\mathbf{1}_S$ for every set $S$ of vertices.
Indeed, we could strengthen this condition by requiring this to hold for all
vectors in $C^0(G;\R)$, not just indicators of sets of vertices. The resulting
relationship between the Laplacians of $G$ and $H$ is described by the Loewner
order on positive semidefinite matrices.
\begin{definition} 
    The \style{Loewner order} on the cone of symmetric positive definite
    matrices of size $n \times n$ is given by the relation $A \preceq B$ if $B
    - A$ is positive semidefinite. Equivalently, $A \preceq B$ if and only if
    $x^TAx \leq x^TBx$ for all $x\in \R^n$.
\end{definition}
By the Courant-Fischer theorem, the relation $A \preceq B$ has important
implications for the eigenvalues and eigenvectors of these matrices. In
particular, $\lambda_{{\max}}(B) \geq \lambda_{{\max}}(A)$ and
$\lambda_{\min}(A) \leq \lambda_{\min}(B)$. If $(1-\epsilon)A \preceq B \preceq
(1+\epsilon)B$, the eigenvalues and eigenvectors of $B$ are constrained to be
close to those of $A$. Thus, it is appropriate to call this relationship a
spectral approximation of $A$ by $B$.

Spielman and Srivastava famously used effective resistances to construct
spectral sparsifiers of graphs~\cite{spielman_graph_2008}. This approach extends
to sheaves on graphs:
just as graph Laplacians can be spectrally approximated by Laplacians of sparse
graphs, so too can sheaf Laplacians be spectrally approximated by Laplacians of
sheaves on sparse complexes. 
\begin{theorem}
\label{thm:sparse}
Let $X$ be a regular cell complex of dimension $d$ and $\Fc$ a cosheaf on $X$
with $\dim C_{d-1}(X;\Fc) = n$. Given $\epsilon > 0$ there exists a subcomplex
$X' \subseteq X$ with the same $(d-1)$-skeleton and $O(\epsilon^{-2}n \log n)$
$d$-cells, together with a cosheaf $\Fc'$ on $X'$ such that
$(1-\epsilon)L_\Fc^{d-1} \preceq L_{\Fc'}^{d-1} \preceq
(1+\epsilon)L_{\Fc}^{d-1}$.
\end{theorem}
\begin{proof}
If $\ker L_{\Fc}^{d-1} = 0$, then an equivalent condition to the conclusion is
\[
  \lambda_{\text{max}}((L_\Fc^{d-1})^{-\frac{1}{2}}L_{\Fc'}^{d-1} (L_\Fc^{d-1})^{-\frac{1}{2}}) \leq 1+\epsilon
\]
and
\[
  \lambda_{\text{min}}((L_\Fc^{d-1})^{-\frac{1}{2}}L_{\Fc'}^{d-1} (L_\Fc^{d-1})^{-\frac{1}{2}}) \geq 1-\epsilon.
\]

If $\ker L_{\Fc}^{d-1}$ is nontrivial, we use the pseudoinverse of
$L_\Fc^{d-1}$ and restrict to the orthogonal complement of the kernel. This
only offers notational difficulties, so in the following we will calculate as
if the kernel were trivial.

Consider the restrictions of $\partial_d$ to each $d$-cell $\sigma$, and note that
\[
  \sum_{\dim(\sigma)=d} (\partial_d|_\sigma)(\partial_d|_\sigma)^*
  =
  L^{d-1}_\Fc
  .
\]
For each $d$-cell $\sigma$, we will choose $\sigma$ to be in $X'$ with probability
\[
  p_\sigma = \min(1,4\epsilon^{-2}\log(n) \tr(R_{\text{eff}}(\sigma)))
  .
\]
If $\sigma$ is chosen to be in $X'$, we choose its extension maps to be
$\Fc_{\sigma\face\bullet}' = \frac{1}{\sqrt{p_\sigma}}\Fc_{\sigma \face
\bullet}$.

Let $A_\sigma$ be independent Bernoulli random variables with $\E[A_\sigma] =
p_\sigma$ and let $X_\sigma = \frac{1}{p_\sigma}A_\sigma
(L_\Fc^{d-1})^{-\frac{1}{2}}
(\partial_d|_\sigma)(\partial_d|_\sigma)^*(L_\Fc^{d-1})^{-\frac{1}{2}}$. Note
that by construction $\sum_\sigma X_\sigma =
(L_\Fc^{d-1})^{-\frac{1}{2}}L_{\Fc'}^{d-1} (L_\Fc^{d-1})^{-\frac{1}{2}}$ and

\[
  \E\left[\sum_\sigma X_\sigma\right]
  =
  (L_\Fc^{d-1})^{-\frac{1}{2}}L_{\Fc}^{d-1} (L_\Fc^{d-1})^{-\frac{1}{2}} = \id_{C_{d-1}(X;\Fc)}
\]

We wish to show that the eigenvalues of $\sum_\sigma X_\sigma$ are close to
those of its expectation, for which we use a matrix Chernoff bound proven
in~\cite{tropp_user-friendly_2012}. This bound requires a bound on the norms of
$X_\sigma$:

\begin{multline*}
\norm{X_\sigma}
\leq
\frac{1}{p_\sigma} \norm{(L_\Fc^{d-1})^{-\frac{1}{2}} (\partial_d|_\sigma)(\partial_d|_\sigma)^*(L_\Fc^{d-1})^{-\frac{1}{2}}} \\
\leq
\frac{1}{p_\sigma} \tr((L_\Fc^{d-1})^{-\frac{1}{2}} (\partial_d|_\sigma)(\partial_d|_\sigma)^*(L_\Fc^{d-1})^{-\frac{1}{2}}) \\
\leq
\frac{1}{p_\sigma} \tr((\partial|_\sigma)^*(L_{\Fc}^{d-1})^\dagger \partial_\sigma)
=
\frac{\tr(R_{\text{eff}}(\sigma))}{p_{\sigma}}.
\end{multline*}

We can {\em a priori} subdivide any $X_\sigma$ with $p_\sigma = 1$ into
sufficiently many independent random variables so that their norms are as small
as necessary. This does not affect the hypotheses of the concentration
inequality, so we consider the case where $p_\sigma < 1$, where
$\norm{X_\sigma} \leq \frac{\epsilon^2}{4 \log n}$. Our matrix Chernoff bound
then gives

\begin{align*}
\Pb\left[\lambda_{\text{min}}((L_\Fc^{d-1})^{-\frac{1}{2}}L_{\Fc'}^{d-1} (L_\Fc^{d-1})^{-\frac{1}{2}}) \leq 1-\epsilon\right] &\leq n \exp\left(\frac{-4\epsilon^2}{2} \epsilon^{-2}\log n\right) = n^{-1}\\
\Pb\left[\lambda_{\text{max}}((L_\Fc^{d-1})^{-\frac{1}{2}}L_{\Fc'}^{d-1} (L_\Fc^{d-1})^{-\frac{1}{2}}) \geq 1+\epsilon \right] &\leq n \exp\left(\frac{-4\epsilon^2}{3} \epsilon^{-2}\log n\right) = n^{-1/3}.
\end{align*}
When $n$ is not trivially small there is therefore a high probability of $L_{\Fc'}$ $\epsilon$-approximating $L_{\Fc}$.

We now check the expected number of $d$-cells in $X'$. This is
\[
  \sum_\sigma p_\sigma \leq 4 \epsilon^{-2}\log n \sum_\sigma \tr(R_{\text{eff}}(\sigma)),
\]
and
\begin{multline*}
\sum_\sigma \tr(R_{\text{eff}}(\sigma))
=
\sum_\sigma \tr((\partial|_\sigma)^*(L_{\Fc}^{d-1})^{-1} \partial|_\sigma)
=
\sum_\sigma \tr((L_{\Fc}^{d-1})^{-1} \partial_\sigma(\partial|_\sigma)^*) \\
=
\tr\left((L_{\Fc}^{d-1})^{-1} \sum_\sigma \partial_\sigma(\partial|_\sigma)^*\right)
=
\tr((L_{\Fc}^{d-1})^{-1} L_\Fc^{d-1}) \leq n.
\end{multline*}

A standard Chernoff bound argument with the Bernoulli random variables
determining whether each cell is included then shows that the number of
$d$-cells is concentrated around its expectation and thus can be chosen to be
$O(\epsilon^{-2}n \log n)$.  
\end{proof}

The proof given here follows the outline of the proof for graphs given by
Spielman in~\cite{spielman_course_2015}. This newer proof simplifies the
original proof in~\cite{spielman_graph_2008}, which used a sampling of edges with replacement. Theorem
\ref{thm:sparse} generalizes a number of theorems on sparsification of graphs
and simplicial complexes
\cite{spielman_spectral_2011,chung_ranking_2012,osting_towards_2017}; however,
it is not the most general sparsification theorem. Indeed, the core argument
does not rely on the cell complex structure, but only on the decomposition of
the Laplacian into a sum of matrices, one corresponding to each cell.

More general, and stronger, theorems about sparsifying sums of symmetric
positive semidefinite matrices have been proven, such as the following
from~\cite{silva_sparse_2016}:

\begin{theorem}[Silva et al. 2016]

Let $B_1,\ldots, B_m$ be symmetric, positive semidefinite matrices of size
$n\times n$ and arbitrary rank. Set $B := \sum_i B_i$. For any $\epsilon \in
(0,1)$, there is a deterministic algorithm to construct a vector $y \in \R^m$
with $O(n/\epsilon^2)$ nonzero entries such that $y$ is nonnegative and
\[
  B \preceq \sum_{i=1}^m y_i B_i \preceq (1+\epsilon)B.
\]
The algorithm runs in $O(mn^3/\epsilon^2)$ time. Moreover, the result continues
to hold if the input matrices $B_1,\ldots, B_m$ are Hermitian and positive
semidefinite.
\end{theorem}

The sheaf theoretic perspective, though not the most general or powerful
possible, nevertheless maintains both a great deal of generality along with a
geometric interpretation of sparsification in terms of an effective resistance.
This geometric interpretation may then be pursued to develop efficient methods
for approximating the effective resistance and hence fast algorithms for
sparsification of cell complexes and cellular sheaves atop them.

% ==========================================================================================
\section{The Cheeger Inequality}
\label{sec:cheeger}
% ==========================================================================================
% -------------------------------------------------------------------------------------
\subsection{The Cheeger Inequality for $O(n)$-bundles}
\label{sec:O(n)bundles}
% -------------------------------------------------------------------------------------

Recall from \S\ref{sec:discvect} the notion of an $O(n)$-bundle on a graph.
Bandeira, Singer, and Spielman proved an analogue to the graph Cheeger
inequality for $O(n)$-bundles~\cite{bandeira_cheeger_2013}. Their goal was to
give guarantees on the performance of a spectral method for finding an
approximate section to a principal $O(n)$-bundle over a graph. For an
$O(n)$-bundle on a graph with Laplacian $L$, degree matrix $D$ and normalized
Laplacian $\mathcal L = D^{-1/2}LD^{-1/2}$, they defined the
\style{frustration} of a $0$-cochain $x$ to be

\[
  \eta(x) = \frac{\ip{x,Lx}}{\ip{x,Dx}} = \frac{\ip{x,\mathcal L x}}{\ip{x,x}},
\]
and showed that any $0$-cochain can be rounded to one with controlled
frustration as follows. Given a $0$-cochain $x$ and a threshold $\Thresh \geq
0$, let $x^\Thresh$ be the cochain whose value at a vertex $v$ is
$x_v/\norm{x_v}$ if $\norm{x_v}^2 \geq \Thresh$ and is zero otherwise. For any
such $x$ there exists a $\Thresh$ such that $\eta(x^\Thresh) \leq \sqrt{10
\eta(x)}$. Taking the minimum over all $0$-cochains $x$ then implies that if
$\lambda_1(\mathcal L)$ is the smallest eigenvalue of $\mathcal L$,

\begin{equation}
\label{eqn:cheeger_bss}
  \lambda_1(\mathcal L) \leq \min_{\norm{x_v} = 1 \text{ or } 0} \eta(x) \leq \sqrt{10 \lambda_1(\mathcal L)}.
\end{equation}

A natural question is whether this theorem extends to a more general class of
sheaves on graphs. One reasonable candidate for extension is the class of
sheaves on graphs where all restriction maps are partial isometries, \ie, maps
which are unitary on the orthogonal complement of their kernels. One might view
these as $O(n)$-bundles where the edge stalks have been reduced in dimension by
an orthogonal projection. However, the cochain rounding approach does not work
for these sheaves, as the following simple counterexample shows. Let $G$ be a
graph with two vertices and one edge, and let $\Fc(v_1) = \Fc(v_2) = \R^2$ and
$\Fc(e) = \R$. Then let $\Fc_{v_1\face e} = [1 \,\, 0]$ and $\Fc_{v_2\face e} =
[\frac{1}{2} \, \, \frac{\sqrt 3}{2}]$. Then let $x_{v_1} = \begin{bmatrix}
    \frac{1}{2} \\ 0\end{bmatrix}$ and $x_{v_2} = \begin{bmatrix} 1 \\ 0
\end{bmatrix}$. Then $\eta(x) = 0$, but $\eta(x^u) > 0$ for any choice of $u <
1$. This means that there cannot exist any function $f: \R \to \R$ with $f(0) =
0$ such that $\eta(x^u) \leq f(\eta(x))$.

This example does not immediately show that the Cheeger
inequality~\eqref{eqn:cheeger_bss} is false for this class of sheaves, since
this sheaf does have a section of stalkwise norm 1, but it does offer a
counterexample to the key lemma in the proof. Indeed, a more complicated family
of counterexamples exists with sheaves that have no global sections. These
counterexamples show that an approach based on variational principles and
rounding is unlikely to prove an analogue of the results of Bandeira et al.\
for more general classes of sheaves.

% -------------------------------------------------------------------------------------
\subsection{Toward a Structural Cheeger Inequality}
\label{sec:structural_cheeger}
% -------------------------------------------------------------------------------------

Many extensions of the graph Cheeger inequality view it from the perspective of
a constrained optimization problem over cochains. This is the origin of the
Cheeger inequality for $O(n)$-bundles, and of the higher-dimensional Cheeger
constants proposed by Gromov, Linial and Meshulam, and
others~\cite{linial_homological_2006,gromov_singularities_2010,parzanchevski_isoperimetric_2016}.
However, a sheaf gives us more structure to work with than simply cochains.

The traditional Cheeger inequality for graphs is frequently stated as a graph
cutting problem: what is the optimal cut balancing the weight of edges removed
with the sizes of the resulting partition? If we take the constant sheaf on a
graph $G$, we can represent a cut of $G$ by setting some restriction maps to
$0$, or, more violently, setting the edge stalks on the cut to be
zero-dimensional. Thus, a potential analogue to the Cheeger constant for
sheaves might be an optimal perturbation to the structure of the sheaf
balancing the size of the perturbation with the size of the support of a new
global section that the perturbation induces.

For instance, we might measure the size of a perturbation of the sheaf's
restriction maps in terms of the square of the Frobenius norm $\norm{\cdot}_F$
of the coboundary matrix.  If we minimize $\norm{\delta_\Fc -
\delta_{\Fc'}}_F^2$, a natural relaxation to the space of all matrices shows
that this value is greater than $\lambda_1(L_\Fc)$.

% ==========================================================================================
\section{Toward Applications}
\label{sec:applications}
% ==========================================================================================

The increase in abstraction and technical overhead implicit in lifting spectral
graph theory to cellular sheaves is nontrivial. However, given the utility of
spectral graph theory in so many areas, the generalization to sheaves would
appear to be a good investment. As this work is an initial survey of the
landscape, we quickly sketch a collection of potential applications. These sketches
are brief enough to allow the curious to peruse, while providing experts with enough
to construct details as needed. 

% -------------------------------------------------------------------------------------
\subsection{Distributed Consensus}
\label{sec:consensus}
% -------------------------------------------------------------------------------------

Graph Laplacians and adjacency matrices play an important role in the study and
design of distributed dynamical systems. This begins with the observation that
the continuous-time dynamical system on a set of real variables 
$\{x_v\colon v\in V(G)\}$ indexed by the vertices of a graph $G$ with Laplacian $L$,
\[
  \dot x = - L x ,
\]
is local with respect to the graph structure: the only terms that influence
$\dot{x}_v$ are the $x_w$ for $w$ adjacent to $v$. Further, if the graph is
connected, diagonalization of $L$ shows that the flow of this dynamical system
converges to a \style{consensus} --- the average of the initial condition. 
A similar observation holds for sheaves of finite-dimensional vector spaces on graphs:

\begin{proposition}

Let $\Fc$ be a sheaf on a cell complex $X$. The dynamical system $\dot x = -
\Delta^k_\Fc x$ has as its space of equilibria $\mathcal{H}^k(X;\Fc)$, the
space of harmonic $k$-cochains of $\Fc$. The trajectory of this dynamical
system initialized at $x_0$ converges exponentially quickly to the orthogonal
projection of $x_0$ onto $\mathcal{H}^k(X;\Fc)$.

\end{proposition}
\begin{proof}

    This is a linear dynamical system with flows given by $x(t) =
    e^{-t\Delta^k_\Fc}x_0$. Since $\Delta^k_\Fc$ is self-adjoint, it has an
    orthogonal eigendecomposition $\Delta^k_\Fc = V\Lambda V^*$, so that flows
    are given by $x(t) = V e^{-t\Lambda}V^*x_0 = \sum_i e^{-t \lambda_i}
    \ip{v_i,x_0} v_i$. The terms of this sum for $\lambda_i > 0$ converge
    exponentially to zero, while the terms with $\lambda_i = 0$ remain
    constant, so that the limit as $t\to \infty$ is $\sum_{v_i \in
    \mathcal{H}^k(X;\Fc)} \ip{v_i,x_0}v_i $. Since the $v_i$ are orthonormal,
    this is an orthogonal projection onto $\mathcal{H}^k(X;\Fc)$. 
\end{proof}

In particular, for $k = 0$, this result implies that a distributed system can
reach consensus on the nearest global section to an initial condition. 

% -------------------------------------------------------------------------------------
\subsection{Flocking}
\label{sec:flocking}
% -------------------------------------------------------------------------------------
One well-known example of consensus on a sheaf 
comes from flocking models~\cite{JLM,tanner_stable_2003}.
In a typical setting, a group of autonomous agents is tasked with arranging themselves 
into a stable formation. As a part of this, agents need a way to agree on a global 
frame of reference.

Suppose one has a collection of autonomous agents in $\R^3$, each of which has
its own internal coordinate system with respect to which it measures the outside
world. Each agent communicates with its neighbors, communication being encoded as a graph. 
Assume agents can calculate bearings to their neighbors in their own coordinate frames. 
The agents wish to agree on a single direction in a global, external frame, perhaps in order 
to travel in the same direction. However, the transformations between local frames are not known.

To solve this, one constructs a sheaf on the neighborhood graph of these agents. The vertices
have as stalks $\R^3$, representing vectors in each agent's individual coordinate frame. 
The edges have stalks $\R^1$, used to compare bearings.  
Since the agents can measure the bearing to each neighbor,
they can project vectors in their coordinate frame onto this bearing. Let
$b(v,w)$ be the unit vector in $v$'s frame pointing toward $w$. Then for the
(oriented) edge $e = v \sim w$, the restriction map $\Fc_{v \face e}: \R^3 \to
\R$ is given by $\ip{b(v,w),\bullet}$, while the restriction map $\Fc_{w \face e}$
is $-\ip{b(w,v),\bullet}$. (The change of sign is necessary because the bearing
vectors are opposites in the global frame.) Any globally consistent direction
for the swarm will be a section of this sheaf, and with a bit more
information the agents can achieve consensus on a direction.

This is merely one simple example of a sheaf for building consensus via Laplacian flow. 
The literature on flocking is quite involved, with various refinements and alternate scenarios,
including, e.g., the case in which the network communication graph changes over time. 

% -------------------------------------------------------------------------------------
\subsection{Opinion Dynamics}
\label{sec:opinion}
% -------------------------------------------------------------------------------------
Sheaf Laplacians provide a drop-in replacement for graph Laplacians when one
wishes to constrain a distributed algorithm to a locally definable subspace of
the global state space. For example, the flocking application in the previous example generalizes
greatly to the setting of opinion dynamics on social networks. 

Consider the setting of a collection of agents, each of whom has an $\R$-valued opinion on some
matter (say, a measure of agreement or disagreement with a particular proposition). A social
network among agents, modeled as a weighted graph, permits influence of opinions based
on Laplacian dynamics: each member of the network continuously adjusts their opinion to be 
more similar to the average of their neighbors' opinions, either in continuous or discrete time. 

The literature on opinion dynamics begins with this simple setting~\cite{IMGroot,Lehrer}, and quickly 
grows to include a number of generalizations of graph Laplacians when there are multiple
opinions or other features~\cite{YTLAA,PB,PiSu}. 
Our perspective is that lifting the simple model to a sheaf automatically incorporates
various novel features while maintaining the simplicity of a Laplacian flow. 

The first obvious generalization is that multiple opinions reside in
higher-dimensional stalks. Each agent $p$ has a vector space $\Fc(p)$ of
opinions on some set of topics. {\em These vector spaces need not have the same
dimension across the network} --- there is no need to assume that every member of
the network has an opinion on every topic. Between each pair $(p,q)$ of
participants adjacent in the social network, there is an edge $e$ with a stalk
$\Fc(e)$, which we might label a \style{discourse space}. Opinions in $\Fc(p)$ and
$\Fc(q)$ are translated into the discourse space by the restriction maps $\Fc_{p
\face e}$ and $\Fc_{q \face e}$. 

The flexibility of a sheaf permits a wonderful array of novel features. For instance, 
each pair of participants in the social network might
only communicate about a handful of topics, and hence only influence each
others' opinions along certain directions. Other topics would therefore lie in
the kernels of the restriction map to their shared discourse space. 

Some features which appear difficult to model in the classical literature on opinion dynamics
are easily programmed into a sheaf: what happens if certain agents lie about their opinion, 
and then only to certain individuals on certain topics? Is a ``public'' consensus (with 
privately-held or context-dependent personal opinions) still possible? This demonstrates the 
utility of sheaves not merely in having stalks which vary, but with varying and interesting restriction maps as well.  

% -------------------------------------------------------------------------------------
\subsection{Distributed Optimization}
\label{sec:opt}
% -------------------------------------------------------------------------------------
Laplacian dynamics are useful not only for mere consensus, but also as a way to
implement consistency constraints for other sorts of distributed algorithms.
Particularly important among these are distributed optimization algorithms,
where a network optimizes a sum of objective functions distributed across the
nodes, with the Laplacian dynamics enforcing the constraint that local state be
consistent across nodes.

Fix a sheaf of finite-dimensional vector spaces over a graph $G$. For each node, $v\in V(G)$, 
assume a cost function, say, a convex functional $\phi_v$ from the stalk of $v$ to $\R$. 
The problem of finding a global section $x=(x_v)$ which minimizes $\sum_v\phi_v(x_v)$ 
subject to the constraint that $x\in H^0(G;\Fc)$ is a relatively unexplored class of
optimization problems. Such problems are naturally distributed in nature, as the constraint
(given in terms of the coboundary operator) is locally defined.   

More generally, one can consider distributed optimization with homological
constraints on any cellular sheaf of vector spaces over a cell complex $X$. 
If the problem to be solved is the optimization of an objective
function defined on $C^k(X;\Fc)$ subject to the constraint that the optimum lie
in $\mathcal{H}^k(X;\Fc)$, then this is naturally distributed. 
One might call such problems \style{homological programs}, analogous to the 
manner in which linear constraints give rise to linear programs. 
The Laplacian evolution then plays a role in providing distributed algorithms to
solve such optimization problems.

% -------------------------------------------------------------------------------------
\subsection{Communication Compression}
\label{sec:comms}
% -------------------------------------------------------------------------------------

Both discrete-time as well as continuous-time Laplacian evolution is useful. 
Consider the following modification of the consensus problems previously discussed. 
Suppose one has a distributed system modeled by a graph $G$,
where each node of $G$ has state in $\R^D$ for some large $D$ and the nodes are
required to reach consensus. The Laplacian flow on states may be discretized as
\begin{equation}
\label{eq:disc-time-lap}
  x[t+1] = (I - \alpha L)x[t].
\end{equation}
To implement this discrete-time evolution equation, at each time step, a node $v$
must send its state $x_v[t]\in \R^D$ to each of its
neighbors, the cost of doing so scaling with state size $D$. It may be preferable
instead to have each node send a lower-dimensional \style{compression} of its state
to each neighbor, \ie, a projection $P_e x_v[t]$ onto a $d\ll D$-dimensional
subspace, with the subspace depending on the edge $e$. 

As with the continuous-time evolution, equilibria of (\ref{eq:disc-time-lap}) correspond
to global sections of the sheaf. However, changing the stalks over edges and the corresponding restriction maps changes the sheaf and therefore, potentially, its global sections. 
The goal for reducing the communication complexity is therefore to program this compressed
sheaf so as to preserve the zeroth cohomology $H^0$. Such a sheaf would comprise a certain {\em approximation} of the constant sheaf over the graph.

% -------------------------------------------------------------------------------------
\subsection{Sheaf Approximation}
\label{sec:approx_sheaf}
% -------------------------------------------------------------------------------------
A sheaf of vector spaces can be thought of as a distributed system of linear
transformations, and its cohomology $H^\bullet$ consists of equivalence classes
of solutions to systems based on these constraints. From this perspective,
questions of approximation --- of sheaves and sheaf cohomology --- take on
especial relevance. The question of approximating {\em global sections} to a given
sheaf has appeared in, \eg,
Robinson~\cite{robinson_sheaves_2017,robinson_assignments_2018}.

Questions of approximating sheaves are equally interesting. Given the relative
lack of investigation, the following definition is perhaps premature;
nevertheless, it is well-motivated by problems of distributed consensus and by 
cognate notions of cellular approximation in algebraic topology.

\begin{definition}
\label{def:sheafappx}

Let $X$ be a regular cell complex, and let $\Gc$ be a sheaf on $X$. We say that
a sheaf $\Fc$ on $X$ is a $k$-\style{approximation} to $\Gc$ if there exists a
sheaf morphism $a: \Gc \to \Fc$ which is an isomorphism on stalks over cells of
degree at most $k$, and which induces an isomorphism $H^i(X;\Gc) \to
H^i(X;\Fc)$ for all $i \leq k$.

\end{definition}

If $\Vspace$ is a vector space, we denote the constant sheaf with stalk
$\Vspace$ by $\constsh{\Vspace}$, and say that $\Fc$ is an \emph{approximation
to the constant sheaf} if $\Fc$ is an approximation to $\constsh{\Vspace}$.

This definition is reminiscent of cellular approximation methods from algebraic
topology. A space $X$ may be $k$-approximated by a cell complex $Y$, via a
morphism $Y \to X$ inducing an isomorphism on homotopy groups up to degree $k$.
Here we approximate a sheaf on a cell complex by one with the same cohomology
in degrees up to $k$.
In particular, a $0$-approximation to $\Fc$ has the same vertex stalks and the
same space of global sections as $\Fc$.

\begin{proposition}\label{prop:approx_const}

If $\Fc$ is a 0-approximation to $\constsh{\Vspace}$, then it is isomorphic to a
sheaf with vertex stalks $\Vspace$ where for each edge $e$ joining vertices $v$
and $w$, the restriction maps $\Fc_{v\face e}: \Vspace \to \Fc(e)$ and
$\Fc_{w\face e}: \Vspace \to \Fc(e)$ are equal.

\end{proposition}
\begin{proof}

Note that because $a: \constsh{\Vspace} \to \Fc$ is an isomorphism on vertex
stalks, $\Fc$ is clearly isomorphic to a sheaf with vertex stalks $\Vspace$.
For every edge $e = (v,w)$ we have the diagram
\[
\begin{tikzcd}
  \Vspace \arrow[d,"\id"'] \arrow[r,"\id"] & \Vspace \arrow[d,"\Fc_{v\face e}"] \\
  \Vspace \arrow[r,"a_e"] & \Fc(e) \\
  \Vspace \arrow[u,"\id"] \arrow[r,"\id"] & \Vspace \arrow[u,"\Fc_{w\face e}"'] \\
\end{tikzcd},
\]
and the only way it can commute is if $\Fc_{v\face e} = \Fc_{w\face e} = a_e$.
\end{proof}

The proof of this proposition shows that specifying an approximation to
$\constsh{\Vspace}$ is the same as specifying a morphism $a_e: \Vspace \to
\Fc(e)$ for each edge $e$ of $G$. Further, in order to produce an approximation
to $\constsh{\Vspace}$, the $a_e$ must assemble to a map $a:
C^1(G;\constsh{\Vspace}) \to C^1(G;\Fc) = \bigoplus_{e \in E} \Fc(e)$ such that
$\ker(a \circ \delta_{\constsh{\Vspace}}) = \ker \delta_{\constsh{\Vspace}}$.
This holds if $\ker a$ is contained in a complement to $\im \delta$;
equivalently, the projection map $\pi: C^1(G;\constsh{\Vspace}) \to
H^1(G;\constsh{\Vspace})$ must be an isomorphism when restricted to $\ker a$.

This suggests a way to construct an approximation to the constant sheaf. Choose
a subspace $K_e$ of $\constsh{\Vspace} (e)$ for each edge $e$ of $G$ and define
$a_e$ to be the projection map $\Vspace \to \Vspace/K_e$. If $\bigoplus_{e \in
E} K_e$ has the same dimension in $H^1(G;\constsh{\Vspace})$ as in
$C^1(G;\constsh{\Vspace})$, then $a = \bigoplus_{e \in E} a_e$ defines the edge
maps giving an approximation to $\constsh{\Vspace}$. (The vertex maps may be
taken to be the identity.)

The question of when a collection $\{K_e\}$ produces
an approximation to the constant sheaf appears quite subtle, and will be the
subject of future work.

This description of approximations to the constant sheaf has ignored the
question of weights --- that is, what inner products to put on the
quotient spaces $\Vspace/K_e$ --- which will be crucial to the spectral behavior of their
Laplacians and hence to the performance of distributed consensus algorithms
based thereon.

To understand this relationship, consider again the discretization of the
Laplacian flow $x[t+1] = (I- \alpha L_\Fc)x[t]$. The matrix $(I- \alpha L_\Fc)$
has an eigenvalue of 1 with eigenspace equal to $H^0(G;\Fc)$, with all other
eigenvalues less than 1. The optimal convergence rate for this consensus
algorithm is obtained at an $\alpha$ keeping the nontrivial eigenvalues as close
to zero as possible. If $\lambda_{\max}$ is the largest eigenvalue of $L_\Fc$
and $\lambda_{\min}$ the smallest nontrivial eigenvalue of $L_\Fc$, this is
obtained at $\alpha = \frac{2}{\lambda_{\max} + \lambda_{\min}}$, for a
nontrivial spectral radius of $r = \frac{\lambda_{\max}-
\lambda_{\min}}{\lambda_{\max} + \lambda_{\min}}$.

The number of steps necessary to reach a level of disagreement of $\epsilon$ is
thus proportional to $\frac{\log \epsilon}{\log(r)}$.
On a $k$-regular graph with $d$-dimensional edge stalks, the total amount of
communication each node must undertake at each step is proportional to $kd$.
Thus the total communication cost per node is proportional to $kd
\frac{\log{\epsilon}}{\log(r)}$. For the constant sheaf $\constsh{\R^D}$,
the total cost per node is $kD \frac{\log{\epsilon}}{\log(R)}$, where $R =
\frac{\lambda_{\max}(G) -
\lambda_{\min}(G)}{\lambda_{\max}(G)-\lambda_{\min}(G)}$ is the corresponding
spectral radius for consensus over $L_{\constsh{\R^D}}$. 
If 
\[\frac{d \log(R)}{D \log(r)} < 1,\]
the total communication cost for consensus using the approximation to the
constant sheaf will be lower than that for the constant sheaf. 
Preliminary investigation suggests it may be possible to construct
approximations to the constant sheaf that achieve this threshold, but more work
is necessary to develop methods for creating spectrally advantageous approximations to the constant
sheaf. These could then be marshaled to improve the speed and efficiency of
distributed algorithms that involve consensus, such as distributed optimization.

% -------------------------------------------------------------------------------------
\subsection{Synchronization}
\label{sec:sync}
% -------------------------------------------------------------------------------------
The concept of \style{synchronization} in the context of problems with data on
graphs is exemplified in work by Singer on the angular alignment
problem~\cite{singer_angular_2011}. The concept was developed further by
Bandeira in his dissertation~\cite{bandeira_convex_2015}.  The general idea is
to recover information about some set of parameters from knowledge of their
pairwise relationships. The general formulation, due to Bandeira, is as follows:
Given a group $G$, a graph $X$, and a function $f_{ij}: G \to \R$ for each edge
$i \sim j$ of $X$, find a function $g: V(X) \to G$ minimizing $\sum_{i \sim j}
f_{ij}(g(i)g(j)^{-1})$.

Often, the functions $f_{ij}$ are chosen such that they have a unique minimum
at a given element $g_{ij} \in G$. One may view this as originating from a
$G$-principal bundle on a graph, with the group elements $g_{ij}$ defining
transition maps. The desired solution is a section of the bundle, determined up
to a self-action of $G$, but may not exist if there is error in the measured
$g_{ij}$. As a result, we seek an error-minimizing solution, where the error is
measured by the functions $f_{ij}$.

Sheaves on graphs offer a broader formulation: synchronization is the problem
of finding a global section, or approximate global section, of an observed
sheaf on a graph. By choosing sheaves valued in a suitable category, we can
recover the group-theoretic formulation. There is something of a gap between
the natural formulation of many synchronization problems and a sheaf valued in
vector spaces; bridging that gap for synchronization over $O(d)$ is the goal
of~\cite{bandeira_cheeger_2013}.

One of the initial motivating problems for the study of synchronization was the
cryo-electron microscopy alignment problem. The goal is to understand the
configuration of a molecule, represented by a density function on $\R^3$.
Cryo-electron microscopy allows one to measure projections of random rotations
of this function onto a fixed two-dimensional plane. One approach to recovery involves
inferring these unknown rotations from pairwise information. After
taking the Fourier transform of the measured two-dimensional distributions,
each pair of distributions will agree on a one-dimensional subspace,
specifically, the invariant axis for the rotation relating the two orientations
of the molecule. 

Suppose the $i^{\rm th}$ measurement is taken from the molecule in an orientation
$\rho_i \in SO(3)$, so that the transformation between the orientation of
measurement $i$ and measurement $j$ is $\rho_{ij} = \rho_j \rho_i^{-1}$. If $x$ is a vector
in the base orientation frame of the molecule, its representation in the frames
for $i$ and $j$ are $\rho_i x$ and $\rho_j x$, respectively. These two vectors
have the same projection onto the invariant subspace of $\rho_{ij}$. Since
this invariant subspace is precisely the subspace on which the relevant two
projections agree, we can check whether two vectors in the frames for
measurements $i$ and $j$ agree by projecting them onto this subspace. 

A single constraint of this form does not ensure equality of vectors in the
different frames. However, sufficiently many generic such constraints will.  
Combining all these pairwise constraints gives us a sheaf  
with the same form as the sheaf of autonomous agents
discussed in \S\ref{sec:flocking}. Note that the pairwise data here obtained is
not in the form of invertible transformations, but weaker constraints. 
Thus, a major motivating problem for synchronization has a natural
expression in the language of cellular sheaves.

The explicitly sheaf-theoretic formulation of the synchronization
problem suggests a different solution approach. If a
synchronization sheaf has a global section --- as is the case when the data are
internally consistent and uncorrupted by noise --- finding that section is
trivial. The traditional approach to synchronization takes these transition
functions as they are, and seeks an approximate section to the sheaf. On the
other hand, we might try to denoise the measured relationships themselves using
the condition of cycle-consistency. That is, given an observed sheaf, find the
nearest sheaf supporting a global section. A structural Cheeger inequality as
discussed in \S\ref{sec:structural_cheeger} would give spectral insights into
this problem. Deeper understanding would come from study of an appropriate
moduli space of cellular sheaves, which is a direction for future work.

% -------------------------------------------------------------------------------------
\subsection{Consistent Clustering}
\label{sec:coclu}
% -------------------------------------------------------------------------------------
As suggested by Gao et al.~\cite{gao_geometry_2016}, the data of a sheaf on a
graph is useful for more than recovering a global section. The problem of
clustering objects where the similarity measure comes from an explicit matching
or transformation gives extra information. If we stipulate that objects within
a cluster should have consistent transformations along cycles, the problem of
clustering becomes the problem of partitioning a graph into subgraphs, each of
which supports the appropriate space of global sections.

Similar ideas arise in~\cite{gao_diffusion_2016}, which considers
correspondences between surfaces produced using the soft Procrustes distance.
These correspondences are maps between probability distributions on the vertex
sets of discretized surfaces. When these surfaces are meshes with varying
numbers of vertices, these maps are not invertible, but by construction they
are represented by doubly stochastic matrices, and the analogue of the inverse
for such a map is simply the transpose of its corresponding matrix. These sorts
of geometric correspondences are natural to consider in the context of
geometric morphometrics, the field devoted to studying and classifying species
based on their geometric properties.

Gao constructs a matrix he calls the \style{graph horizontal Laplacian}, together with
normalized versions he uses to formulate a fiber bundle version of the
diffusion maps algorithm for dimensionality reduction. The graph horizontal
Laplacian is related to a map of graphs $X \to G$, where the fibers over
vertices of $G$ are discrete. A weighting on the edges of $X$ induces a
matrix-valued weighting on the edges of $G$. This produces a weighted adjacency
matrix $W$ of $G$, from which the graph horizontal Laplacian is generated by
$L^H = D - W$, where $D$ is the diagonal matrix necessary to make $L^H$ have
row sums equal to zero. 

This is in fact equivalent to the sheaf Laplacian of
the pushforward of the weighted constant sheaf on $X$, and as a consequence of
Proposition \ref{prop:pushforward}, is simply a block subdivision of the
Laplacian of $X$. This sheaf on $G$ can then be normalized to construct a
diffusion map embedding of the vertices of $G$, as well as an embedding of the
vertices of $X$. When applied to the surface correspondence problem, the
eigenvectors of the resulting sheaf Laplacian serve to partition the surfaces
into automatically determined landmarks or regions of interest.

Approaching these notions of partitioning, partial sections, and noninvertible
matchings from a sheaf-theoretic perspective offers new tools and clarifies the
problems in question, with potential for the spectral approach to yield insights. 

% ==========================================================================================
\section{Closing Questions}
\label{sec:further}
% ==========================================================================================

There are numerous interesting open questions in an emerging spectral sheaf theory. We
highlight a few below, with comments.

% -------------------------------------------------------------------------------------
\subsection{Metrics on the Space of Cellular Sheaves}
% -------------------------------------------------------------------------------------
Interleaving-type constructions have been used to define metrics on the space
of constructible sheaves. However, these rely on explicit geometric information
about the sheaves and their underlying spaces. Working with weighted cellular
sheaves may make it possible to define useful distances that rely only on the
combinatorial and algebraic structure of the sheaves. What are the most useful
metrics on the space of sheaves? How do they interact with the sheaf Laplacians
and their spectra? How does this shed light on a moduli space of cellular sheaves?

% -------------------------------------------------------------------------------------
\subsection{Developing a Cheeger Inequality}
% -------------------------------------------------------------------------------------
Following the discussion in \S\ref{sec:structural_cheeger}, a structural
Cheeger inequality for sheaves is connected to questions of a potential moduli space of
sheaves. Such an inequality would describe how the spectral properties of a
sheaf interact with its distance to the nearest sheaf with a nontrivial global
section. Can a Cheeger inequality emerge from approximations to the constant
sheaf, seeking $0$-cochains of small coboundary which are constant on large
sets of vertices?

% -------------------------------------------------------------------------------------
\subsection{Interactions with the Derived Category}
% -------------------------------------------------------------------------------------
The standard way of understanding sheaf cohomology is through the derived
category of complexes of sheaves~\cite{gelfand_methods_2003}. We may replace a
sheaf by an injective resolution and take the cohomology of the complex of
sheaves. What is the relationship between a weighted sheaf and its injective
resolutions, and how do the resulting Laplacians connect with the Hodge
Laplacian defined on the cochain complex? What results can be proven about
their spectra? How do they interact with the standard sheaf operations? Is
there a consistent way to add weights to the derived category of sheaves? We
should not expect the answers to these questions to be unique due to the dagger
categorical issues discussed in \S\ref{sec:weighted}, but there may be
constructions which are nevertheless appealing.

% -------------------------------------------------------------------------------------
\subsection{Random Walks}
% -------------------------------------------------------------------------------------
Chung and Zhao~\cite{chung_ranking_2012} considered random walks on discrete
$O(n)$-bundles, including a definition of a sort of \style{PageRank} algorithm.
Is it possible to define random walks on more general sheaves of vector spaces, 
and to what extent are such related to the sheaf Laplacian and its spectral features? 
Is there an analogous PageRank algorithm for sheaves?

% -------------------------------------------------------------------------------------
\subsection{Cones and Directedness}
% -------------------------------------------------------------------------------------
How does one model directedness and asymmetric relations on sheaves? Sheaves of
cones and sheaf cohomology taking values in categories of cones have proven
useful in recent applications of sheaf theory to problems incorporating
directedness~\cite{ghrist_positive_2017,kashiwara_persistent_2018}. Such
methods, though promising, may be noncommutative, using semigroups or
semimodules to encode the directedness, which, in turn, pushes the boundaries
of existing methods in sheaf theory and nonabelian sheaf cohomology.

\bibliographystyle{alpha}
%\bibliography{sheafspectra}

\end{document}